\documentclass{siamart}
\pdfoutput=1

\usepackage{graphicx}
\usepackage{ifthen}
\newboolean{epsfigs}
\setboolean{epsfigs}{false}

\usepackage{color}

\def\mydate{\number\day\ {\ifcase\month \or January\or February\or
              March\or April\or May\or June\or July\or August\or
              September\or October\or November\or December\fi}
\number\year}

\ifthenelse{\boolean{epsfigs}}{\input{SETEPS}}{}

\usepackage{bbm}
\def\C{\mathbbm{C}}

\def\R{\mathbbm{R}}
\def\Rn{\R^n}

\def\W{W}      
\def\wh#1{\widehat{#1}}

\def\poly{\varphi}
\def\mingmres{\min_{\stackrel{\scriptstyle{\poly\in {\cal P}_k}}
                            {\scriptstyle{\poly(0)=1}}}}
\def\Orsirr{Orsirr\textunderscore 1}
\def\loc{\mbox{\emph{loc}}}
\def\ip#1{\langle #1\rangle}

\newsiamremark{example}{Example}
\mathcode`\@="8000
{\catcode`\@=\active \gdef@{\mkern1mu}} 

\begin{document}
\title{Weighted Inner Products for GMRES and GMRES-DR}

\author{
Mark Embree\footnotemark[1]
\and Ronald B. Morgan\footnotemark[2]
\and Huy V. Nguyen\footnotemark[3] }

\maketitle

\renewcommand{\thefootnote}{\fnsymbol{footnote}}
\footnotetext[1]{Department of Mathematics and
Computational Modeling and Data Analytics Division,
Academy of Integrated Science,
Virginia Tech, Blacksburg, VA 24061
(\email {embree@vt.edu}).
Supported in part by National Science Foundation grant DGE-1545362.}
\footnotetext[2]{Department of Mathematics, Baylor
University, Waco, TX 76798-7328\\ (\email{Ronald\_Morgan@baylor.edu}).
Supported in part by the National Science Foundation Computational Mathematics Program 
under grant DMS-1418677.}
\footnotetext[3]{Department of Mathematics, Computer Science, and Cooperative Engineering,
University of St.~Thomas, 3800 Montrose, Houston, TX 77006 (\email{nguyenhv2@stthom.edu}).}
\renewcommand{\thefootnote}{\arabic{footnote}}


\begin{abstract}
The convergence of the restarted GMRES method can be significantly improved, 
for some problems, by using a weighted inner product that changes at each restart.
How does this weighting affect convergence, and when is it useful?
We show that weighted inner products can help in two distinct
ways: when the coefficient matrix has localized eigenvectors, weighting
can allow restarted GMRES to focus on eigenvalues that otherwise slow
convergence; for general problems, weighting can break the cyclic 
convergence pattern into which restarted GMRES often settles.
The eigenvectors of matrices derived from differential equations are 
often not localized, thus limiting the impact of weighting.
For such problems, incorporating the discrete cosine transform 
into the inner product can significantly improve GMRES convergence,
giving a method we call W-GMRES-DCT.\ \ 
Integrating weighting with eigenvalue deflation via GMRES-DR 
also can give effective solutions.
\end{abstract}

\begin{keywords} 
linear equations, localized eigenvectors, weighted inner product, restarted GMRES, GMRES-DR, W-GMRES-DCT, deflation
\end{keywords}

\begin{AMS}
65F10, 15A06
\end{AMS}


\headers{Weighted GMRES and GMRES-DR}{Mark Embree, Ronald B. Morgan, and Huy Nguyen}

\section{Introduction}
We seek to solve large nonsymmetric systems of linear equations 
using variants of the GMRES method~\cite{SaSc} that change the 
inner product each time the method is restarted. 

To solve $Ax = b$ for the unknown $x$, GMRES computes the approximate 
solution $x_k$ from the Krylov subspace 
${\cal K}_k(A,b) = {\rm span}\{b, Ab, \ldots, A^{k-1}b\}$
that minimizes the norm of the residual $r_k := b - A x_k$.
This residual is usually optimized in the (Euclidean) 2-norm, 
but the algorithm can be implemented in any norm induced by an
inner product.  This generality has been studied in 
theory (e.g.,~\cite{FLT08, FM84}) and in practice.
The inner product does not affect the spectrum of $A$, nor the 
subspace ${\cal K}_k(A,b)$ (assuming the method is not restarted),
but it can significantly alter the departure of $A$ from normality and
the conditioning of the eigenvalues.
Numerous authors have proposed inner products in which 
certain preconditioned saddle point systems are self-adjoint, 
enabling optimal short-recurrence Krylov methods in place of the 
long-recurrences required for the 2-norm~\cite{BP88,FRSW98,Kla98,LP08,SW08}.
Recently Pestana and Wathen used non-standard 
inner products to inform preconditioner design~\cite{PW13}.

For any inner product $\ip{\cdot, \cdot}_W$ on $\R^n$, 
there exists a matrix $W\in\R^{n\times n}$ for which $\ip{u,v}_W = v^* W u$ 
for all $u, v\in\Rn$,
and $W$ is Hermitian positive definite in the Euclidean inner product.
(Throughout, $^*$ denotes the conjugate transpose.)
With the $W$-inner product is associated the norm 
$\|u\|_W = \sqrt{\ip{u,u}_W} = \sqrt{u^* W u}$.\ \ 
To implement GMRES in the $W$-inner product, simply
replace all norms and inner products in the standard GMRES algorithm.
Any advantage gained by this 
new geometry must be balanced against the extra expense
of computing the inner products: essentially one
matrix-vector product with $W$ for each inner product and norm evaluation.

In 1998 Essai proposed an intriguing way to utilize
the inner product in restarted GMRES~\cite{Es98}:
take $W$ to be diagonal 
and \emph{change the inner product at each restart}.
His computational experiments show that this 
``weighting'' can significantly improve the performance of
restarted GMRES, especially for difficult problems and frequent restarts.  
Subsequent work has investigated how weighted GMRES performs 
on large examples~\cite{CaYu,GP14,Me09,NiLuZh,SNZa}.
We offer a new perspective on why weighting helps (or not), 
justified by simple analysis and careful computational experiments.

\Cref{sec:wgmres} describes Essai's algorithm.
In \cref{sec:analysis}, we provide some basic analysis 
for weighted GMRES and illustrative experiments, arguing that the 
algorithm performs well when the eigenvectors are \emph{localized}.
In \cref{sec:fft} we propose a new variant that combines weighting 
with the discrete cosine transform, which localizes eigenvectors in
many common situations.\ \ 
\Cref{sec:drwgmres} describes a variant with deflated restarting (GMRES-DR~\cite{GMRES-DR}) 
that both solves linear equations and computes eigenvectors.  
The eigenvectors improve the convergence of the solver for the linear equations, 
and a weighted inner product can help further.

\section{Weighted GMRES} \label{sec:wgmres}
The standard restarted GMRES algorithm~\cite{SaSc} computes 
the residual-minimizing approximation from a degree-$m$ Krylov subspace, 
then uses the residual associated with this solution estimate 
to generate a new Krylov subspace.  We call this algorithm \emph{GMRES($m$)}, 
and every set of $m$ iterations after a restart a \emph{cycle}.
Essai's variant~\cite{Es98}, which we call \emph{W-GMRES($m$)}, 
changes the inner product at the beginning of each cycle of restarted GMRES, 
with weights determined by the residual vector obtained from the last cycle.  
Let $r_{km} = b - A x_{km}$ denote the GMRES($m$) residual after $k$~cycles.
Then the weight matrix defining the inner product for the next cycle is 
$W = {\rm diag}(w_1,\ldots, w_n)$, where
\begin{equation} \label{eq:essai}
  w_j = \max \bigg\{ {|(r_{km})_j| \over \|r_{km}\|_\infty}, 10^{-10}\bigg\},
\end{equation}
where $(r_{km})_j$ denotes the $j$th entry of $r_{km}$.%
\footnote{\label{fn:wsmall}
Essai scales these entries differently, but, as Cao and Yu~\cite{CaYu} note,
scaling the inner product by a constant does not affect the residuals 
produced by the algorithm.}
The weighted inner product is then $\ip{u,v}_W = v^* W u$ 
giving the norm $\|u\|_W = \sqrt{\ip{u,u}_W}$.  
At the end of a cycle, the residual vector $r_{km}$ is computed 
and from it the inner product for the next cycle is built.
We have found it helpful to restrict the diagonal entries in~(\ref{eq:essai}) 
to be at least $10^{-10}$, thus limiting the condition number of $W$ 
to $\|W\|_2 \|W^{-1}\|_2 \le 10^{10}$.
Finite precision arithmetic influences GMRES in nontrivial ways;
see, e.g., \cite[sect.~5.7]{LS13}.
Weighted inner products can further complicate this behavior;
a careful study of this phenomenon is beyond the scope 
of our study, but note that Rozlo\v{z}n\'{\i}k et al.~\cite{RTSK12} 
have studied orthogonalization algorithms in general inner products.
(Although reorthogonalization is generally not needed for GMRES~\cite{PRS06}, 
all tests here use one step of full reothogonalization for extra stability.
Running the experiments in \cref{sec:wgmres} without this step yields
qualitatively similar results, but the iteration counts differ.
All computations were performed using MATLAB 2015b.  Different
MATLAB versions yield somewhat different results, especially
for experiments that require many iterations.)

For diagonal $W$, the inner product takes $3n$ operations 
(additions and multiplications), instead of $2n$ for the Euclidean case. 
This increase of roughly $25\%$ in the Gram--Schmidt step in GMRES 
can be significant for very sparse $A$, but if the matrix-vector product
is expensive, the cost of weighting is negligible.

How can weighting help?  Suppose we restart GMRES every $m$ iterations.  
GMRES minimizes a residual of the
form $\poly(A)r$, where $\poly$ is a polynomial of degree no greater
than $m$ with $\poly(0)=1$.  
When $A$ is nearly normal, one learns much from the 
magnitude of $\poly$ on the spectrum of~$A$.\ \
Loosely speaking, since GMRES must minimize the norm of the 
residual, the optimal $\poly$ cannot target isolated
eigenvalues near the origin, since $\poly$ would then be large at other eigenvalues.
This constraint can cause uniformly slow convergence 
from cycle to cycle; the restarted algorithm
fails to match the ``superlinear'' convergence 
enabled by higher degree polynomials in full GMRES~\cite{vdVV93}.
(Higher degree polynomials can target a few critical 
eigenvalues with a few roots, while having sufficiently many roots 
remaining to be uniformly small on the rest of the spectrum.)
Weighting can skew the geometry 
to favor the troublesome eigenvalues that impede convergence,
leaving easy components to be eliminated quickly in subsequent cycles.
In aggregate, weighting can produce a string of iterates, each sub-optimal 
in the 2-norm, that collectively give faster convergence than 
standard restarted GMRES, with its locally optimal (but globally 
sub-optimal) iterates.

Essai's strategy for the inner product relies on 
the intuition that one should emphasize those 
components of the residual vector that have the largest magnitude.   
Intuitively speaking, these components have thus far been 
neglected by the method, and may benefit from some preferential treatment.
By analyzing the spectral properties of the weighted GMRES algorithm
in \cref{sec:analysis}, we explain when this intuition 
is valid, show how it can go wrong, and offer a possible remedy.

Essai showed that his weighting can give significantly 
faster convergence for some problems, but did not provide 
guidance on matrix properties for which this was the case.
Later, Saberi Najafi and Zareamoghaddam~\cite{SNZa} 
showed that diagonal $W$ with entries
from a \emph{random} uniform distribution on $[0.5,1.5]$ can sometimes be 
effective, too.
Weighting appears to be particularly useful when the 
GMRES($m$) restart parameter $m$ is small, though it can help for
larger $m$ when the problem is difficult.

\begin{figure}[t!]
           \begin{center}\includegraphics[scale=0.55]{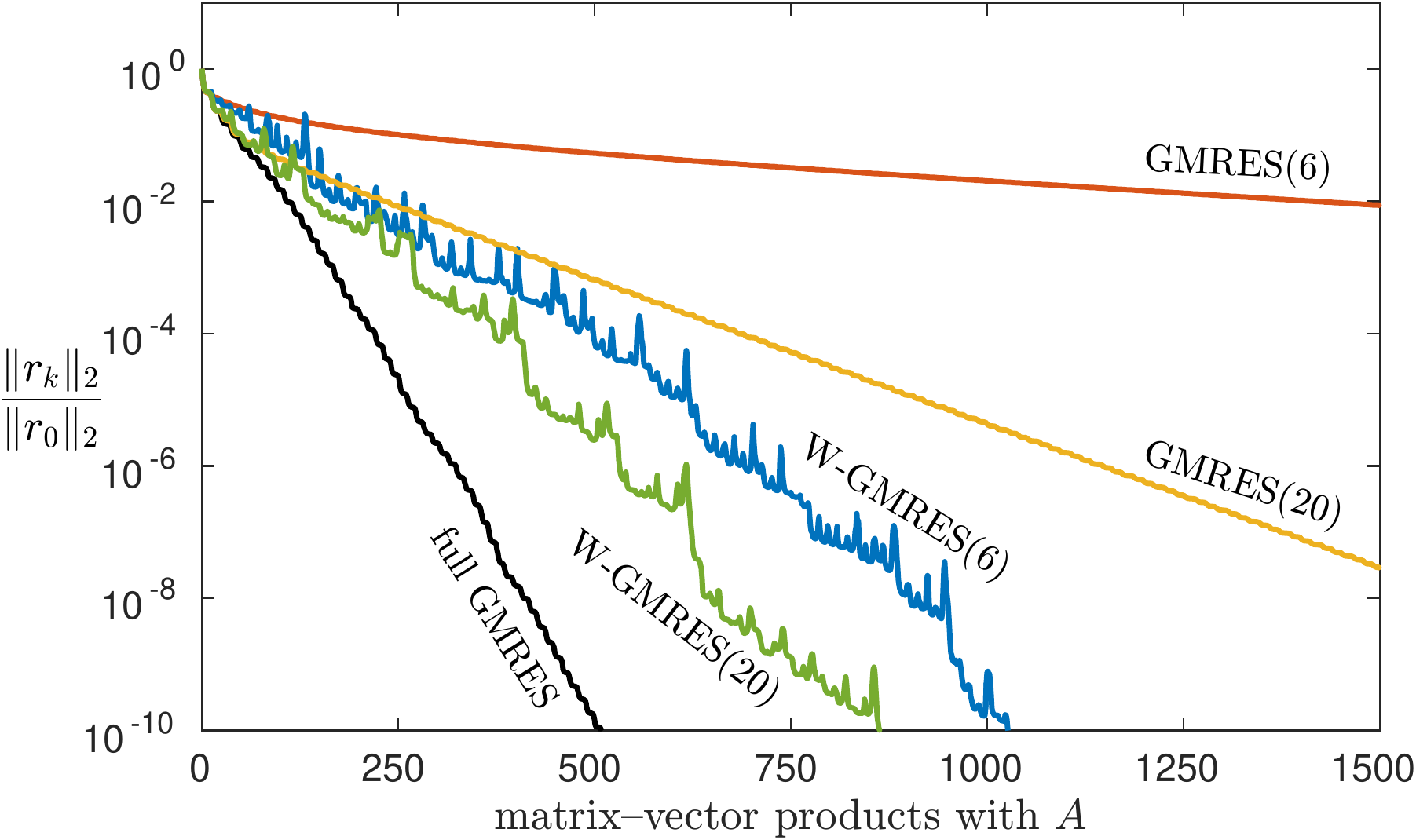} 
            \end{center}

\vspace*{-5pt}
\caption{\label{fig:gwf1}
Convergence of full GMRES, GMRES($m$), and W-GMRES($m$) 
with $m=6$ and $m=20$ for the Add20 matrix.
In this and all subsequent illustrations, the residual norm is measured
in the standard 2-norm for all algorithms.
Full GMRES converges to this tolerance in about 509~iterations.}
\end{figure}

\begin{example}  \label{ex:add20}
Consider the Add20 example from the Matrix Market collection~\cite{MM},
a nonsymmetric matrix of dimension $n=2395$ used by Essai~\cite{Es98}.
(Apparently our ``cycles'' correspond to Essai's ``iterations.'')
Here $b$ is a random Normal(0,1) vector.  
\Cref{fig:gwf1} compares GMRES($m$) and W-GMRES($m$) for $m=6$ and $20$.
In terms of matrix-vector products with $A$, 
W-GMRES(6) converges about~12.5 times faster than GMRES(6);
the relative improvement for $m=20$ is less extreme, but
W-GMRES is still better by more than a factor of two.
(Full GMRES converges to $\|r_k\|_2/\|r_0\|_2 \le 10^{-10}$ in 509~iterations.)
Note that $\|r_k\|_2$ does not converge monotonically for W-GMRES, 
since the residuals are minimized in the weighted norm, not the 2-norm.
\end{example}

Is weighting typically this effective?
Essai shows several other tests where weighting helps, 
though not as much as for Add20.
Cao and Yu~\cite{CaYu} give evidence that suggests
weighting is not as effective for preconditioned problems,
though it can still help.  
G\"uttel and Pestana arrive at a similar conclusion after 
conducting extensive numerical experiments on 109 ILU-preconditioned 
test matrices:
``weighted GMRES may outperform unweighted GMRES for some problems, 
but more often this method is not competitive with other Krylov subspace 
methods\ldots''~\cite[p.~733]{GP14}
and 
``we believe that WGMRES should not be used in combination with 
preconditioners, although we are aware that for some examples it 
may perform satisfactorily''~\cite[p.~750]{GP14}.
 Yet weighting is fairly inexpensive and can improve convergence for some difficult problems.  When weighting works, what makes it work?

\section{Analysis of Weighted GMRES} \label{sec:analysis}

Analyzing W-GMRES($m$) must be no easier than analyzing
standard restarted GMRES, an algorithm for which results
are quite limited.  
(One knows conditions under which GMRES($m$) converges
for all initial residuals~\cite{FJKM96}, and about 
cyclic behavior for GMRES($n-1$) for normal $A$~\cite{BaJeMa}.
Few other rigorous results are known.)
Restarted GMRES is a nonlinear dynamical system 
involving the entries of $A$ and $b$, and
experiments suggest its convergence can 
depend sensitively on $b$~\cite{Em03}.
Aware of these challenges, we seek some insight 
by first studying a setting for which W-GMRES is ideally suited.

Let $\W = S^*S$ denote a positive definite matrix that induces the inner product
$\ip{u,v}_\W = v^* \W u$ on $\Rn$, with induced norm 
\[ \|u\|_\W = \sqrt{\ip{u,u}_\W} = \sqrt{u^*\W u} = \sqrt{u^*S^*Su} = \|S@u\|_2.\]
Given the initial guess $x_0$ and residual $r_0 = b-A x_0$, at the $k$th step GMRES in the 
$W$-inner product computes the residual $r_k$ that satisfies
\begin{eqnarray*}
  \|r_k\|_\W = \mingmres \|\poly(A) r_0\|_\W 
            &=& \mingmres \|S@@ \poly(A)S^{-1} Sr_0\|_2 \\[0.25em]
            &=& \mingmres \|\poly(SAS^{-1}) Sr_0\|_2.
\end{eqnarray*}
(Here ${\cal P}_k$ denotes the set of polynomials of degree $k$ or less.)
Thus, before any restarts are performed, 
 W-GMRES applied to $(A,r_0)$ is equivalent to 2-norm GMRES
applied to $(SAS^{-1}, S@r_0)$.  
(From this observation G\"uttel and Pestana propose an alternative implementation of W-GMRES, 
their Algorithm~2~\cite[p.~744]{GP14}.)
This connection between W-GMRES and standard GMRES 
was first noted in an abstract by Gutknecht and Loher~\cite{GL01}, 
who considered ``preconditioning by similarity transformation.'' 
Since $A$ and $SAS^{-1}$ have a common spectrum, one expects full GMRES 
in the 2-norm and $\W$-norm to have similar \emph{asymptotic} 
convergence behavior.
However, the choice of $S$ can significantly affect the 
departure of $A$ from normality.  
For example, if $A$ is diagonalizable with $A = V \Lambda V^{-1}$,
then taking $\W = V^{-*}V^{-1}$ with $S=V^{-1}$ renders
$SAS^{-1} = \Lambda$ a diagonal (and hence normal) matrix.

The $W$-norms of the residuals produced by a cycle of W-GMRES 
will decrease monotonically, but this need
not be true of the 2-norms, as seen in \cref{fig:gwf1}. 
In fact, 
\[ \|r_k\|_\W = \|S r_k\|_2 \le \|S\|_2 \|r_k\|_2,\]
and similarly
\[ \|r_k \|_2 = \|S^{-1}S r_k\|_2 
     \le \|S^{-1}\|_2 \|S r_k\|_2
     = \|S^{-1}\|_2 \|r_k\|_\W, \]
so 
\[ {1\over \|S\|_2} \|r_k\|_\W \le \|r_k\|_2 \le \| S^{-1}\|_2 @ \|r_k\|_\W.\]
Pestana and Wathen~\cite[Thm.~4]{PW13} show that
the same bound holds when 
the residual in the middle of this inequality is replaced by the optimal 
2-norm GMRES residual.  
That is, if $r_k^{(\W)}$ and $r_k^{(2)}$ denote the $k$th residuals
from GMRES in the $\W$-norm and 2-norm
(both before any restart is performed, $k\le m$), then
\[ {1\over \|S\|_2} \|r_k^{(\W)}\|_\W 
   \le \|r_k^{(2)}\|_2 
\le \| S^{-1}\|_2 @ \|r_k^{(\W)}\|_\W.\]
Thus for $\W$-norm GMRES to depart significantly from 2-norm GMRES
\emph{over the course of one cycle}, $S$ must be (at least somewhat) 
ill-conditioned.
However, restarting with a new inner product can change the
dynamics across cycles significantly.

\subsection{The ideal setting for W-GMRES} \label{sec:diag}

For diagonal $W$, as in~\cite{Es98}, let
\[ S = {\rm diag}(s_1, \ldots, s_n):= W^{1/2}.\]
When the coefficient matrix $A$ is diagonal,
\[ A = {\rm diag}(\lambda_1, \ldots, \lambda_n),\]
Essai's weighting is well motivated and its effects can be readily understood.%
\footnote{Examples~1 and~2 of G\"uttel and Pestana~\cite{GP14} 
use diagonal $A$, and the analysis here can help explain their experimental results.}
In this case, the eigenvectors are columns of the identity matrix,
$V = I$, and $SA=AS$.  
Consider a cycle of W-GMRES($m$) starting with
$r_0 = b = [b_1, \ldots, b_n]^T$.
For $k\le m,$
\begin{subequations} \label{eq:diag}
\begin{eqnarray}
 \|r_k\|_\W^2 = \mingmres \|\poly(A)b\|_\W^2 
            &=& \mingmres \|S@\poly(A) b\|_2^2 \\[.25em]
            &=& \mingmres \|\poly(SAS^{-1})Sb\|_2^2 \\[-.5em]
            &=& \mingmres \|\poly(A)Sb\|_2^2 
            = \mingmres \sum_{j=1}^n  |@\poly(\lambda_j)|^2 |s_j b_j|^2. \label{eq:diag3}\\ \nonumber
\end{eqnarray}
\end{subequations}
This last formula reveals that the weight $s_j$ affects GMRES 
in the same way as the component $b_j$ of the right-hand side
(which is the component of $b$ in the $j$th eigenvector direction for this $A$), 
\emph{so the weights tune how much GMRES will emphasize any given eigenvalue}.
Up to scaling, Essai's proposal amounts to $s_j = \sqrt{|b_j|}$, so 
\begin{equation} \label{eq:diagessai}
  \|r_k\|_\W^2 = \mingmres \sum_{j=1}^n  |@\poly(\lambda_j)|^2 @ |b_j|^3,
\end{equation}
This approach favors large components in $b$ over smaller ones.
Small entries in $b$ (reduced at a previous 
cycle) were ``easy'' for GMRES to reduce because, 
for example, they correspond to eigenvalues far from the origin.  
Large entries in $b$ have been neglected by previous cycles, 
most likely because they correspond to eigenvalues near the origin.
A low-degree polynomial that is small at these points is likely to 
be large at large magnitude eigenvalues, and so will not be picked as 
the optimal polynomial by restarted GMRES.
Weighting tips the scales to favor of these neglected eigenvalues
near the origin for one cycle; this preference will usually 
\emph{increase} the residual in ``easy'' components, 
but such an increase can be quickly remedied at the next cycle.  
Consider how standard GMRES(1) handles the scenario
\[ 0< \lambda_1 \ll 1 \le \lambda_2 \approx \cdots \approx \lambda_n.\]
Putting the root $\zeta$ of the residual polynomial $\poly(z) = 1-z/\zeta$ 
near the cluster of eigenvalues $\lambda_2\approx\cdots\approx \lambda_n$
will significantly diminish the residual vector in components $2, \ldots, n$,
while not increasing the first component, since $|\poly(\lambda_1)|<1$.   
On the other hand, placing $\zeta$ near $\lambda_1$ would give
$|\poly(\lambda_j)|\gg 1$ for $j=2,\ldots, n$, increasing 
the overall residual norm, but in a manner that could be remedied at 
the next cycle if GMRES(1) could take a locally
suboptimal cycle to accelerate the overall convergence.
Essai's weighting enables GMRES to target 
that small magnitude eigenvalue for one cycle, potentially
at an increase to the 2-norm of the residual that can be 
corrected at the next step without undoing the reduced component 
associated with $\lambda_1$.

\begin{figure}[b!]
 \begin{center}\includegraphics[scale=0.5]{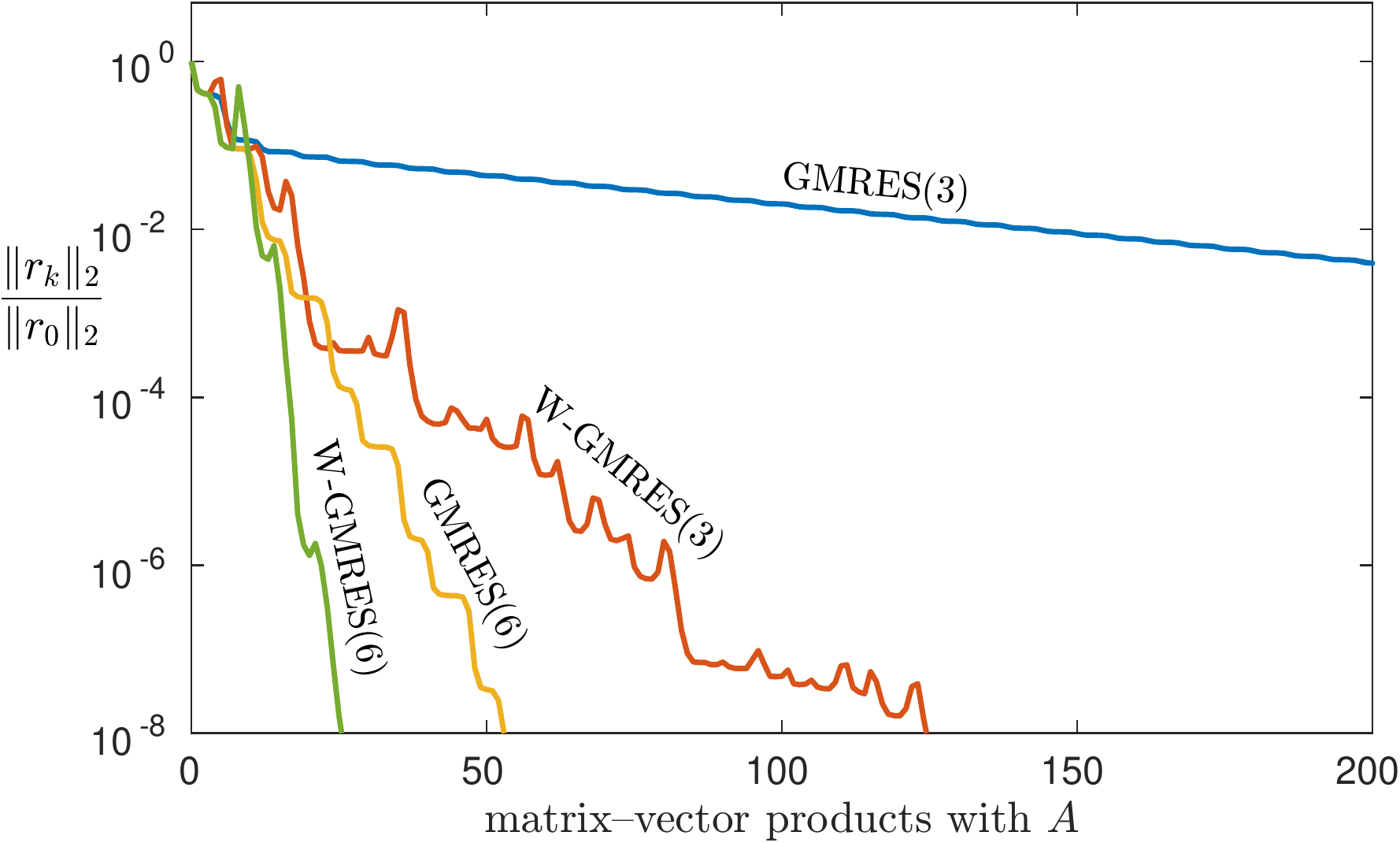} \end{center}

\vspace{-5pt}
\caption{\label{fig:gwf2a}
GMRES($m$) and W-GMRES($m$) with $m=3, 6$ 
for the diagonal test case, \cref{ex:10diag}.}
\end{figure}

\begin{figure}[t!]
\begin{center}
\includegraphics[scale=0.33]{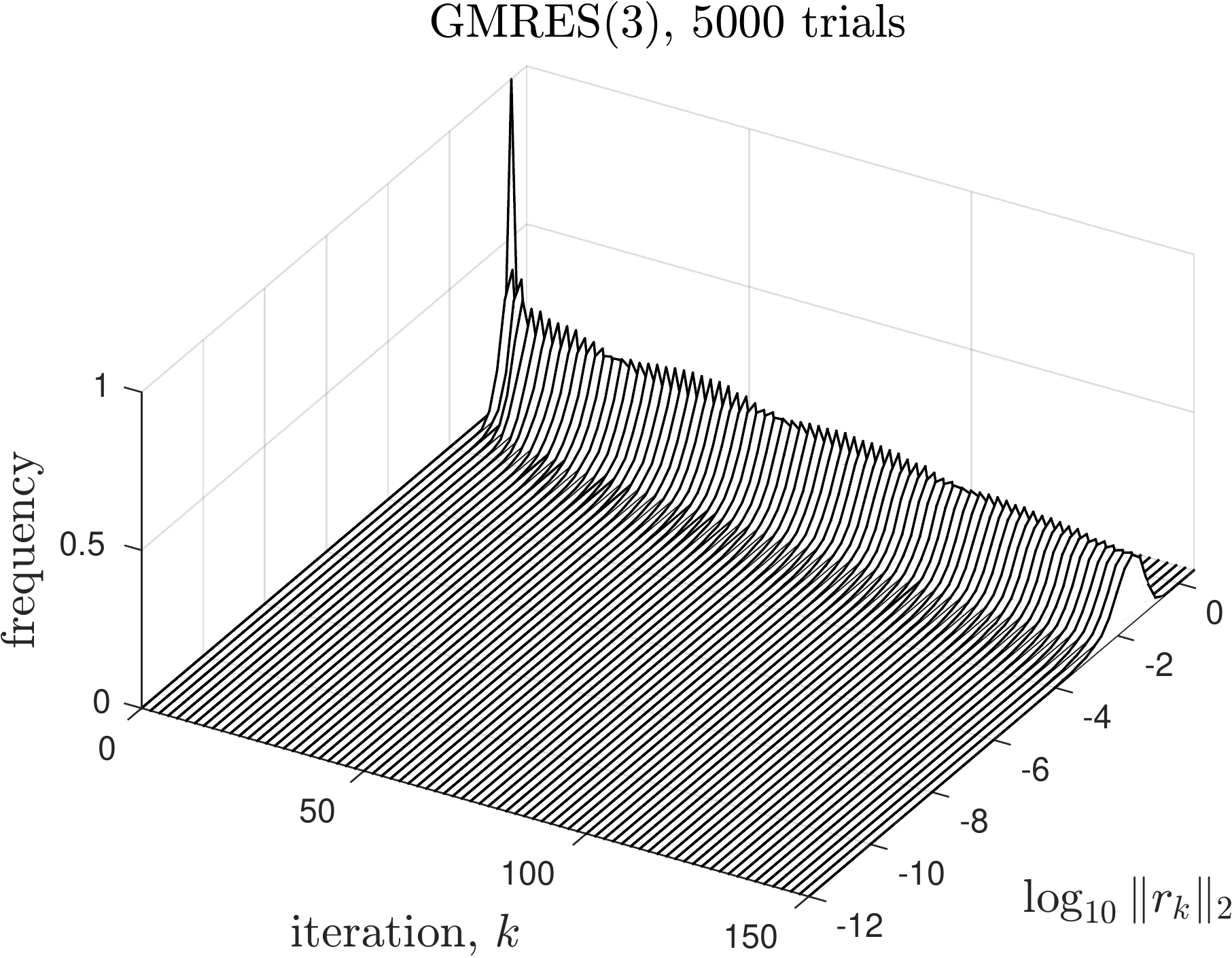}
\qquad
\includegraphics[scale=0.33]{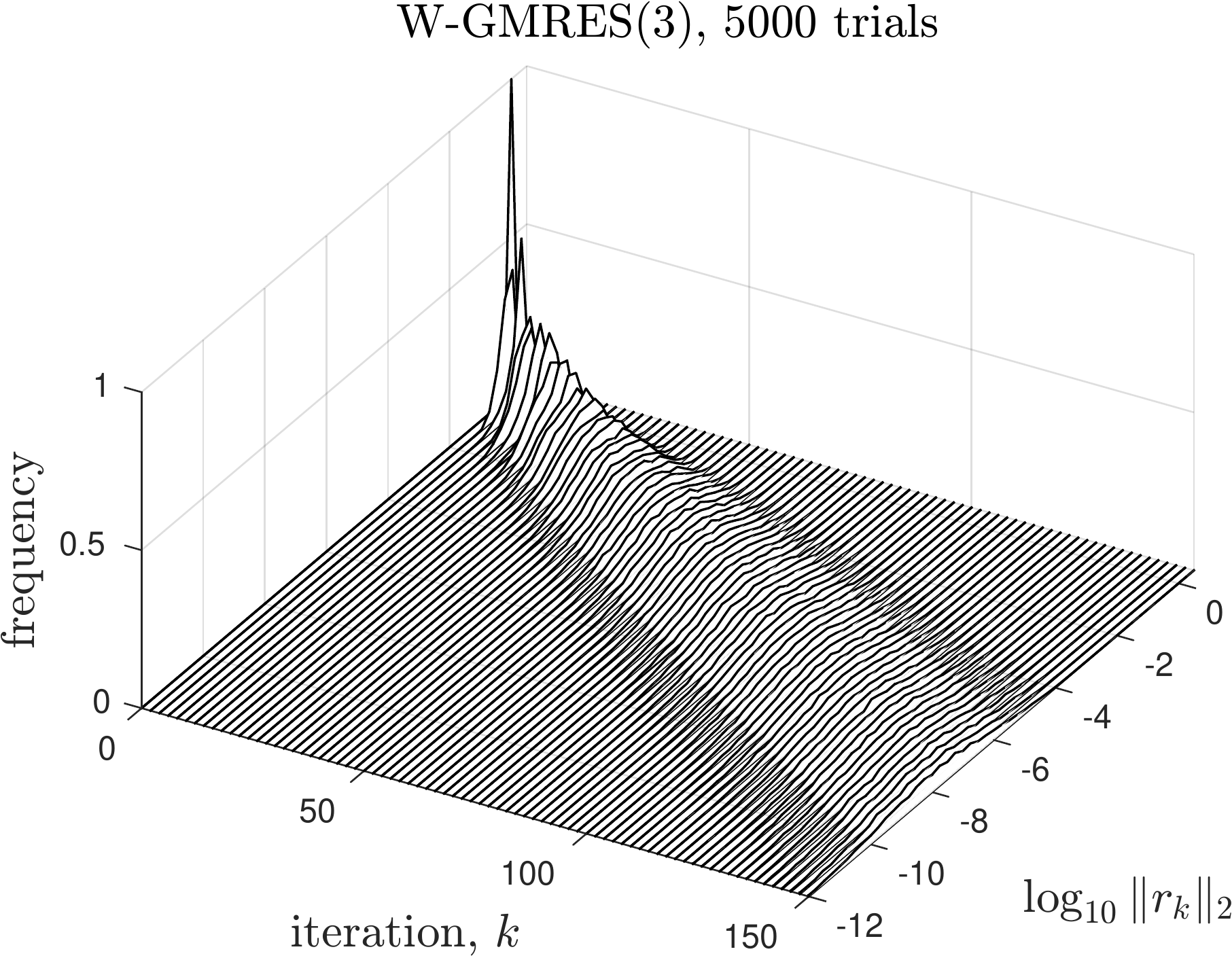}
\end{center}
\caption{\label{fig:hist}
Histograms of GMRES(3) and W-GMRES(3) convergence for 
\cref{ex:10diag} with 5000 random unit vectors $b=r_0$
(the same for both methods).  For each $k$, these plots show the distribution 
of $\|r_k\|_2$; while W-GMRES(3) converges faster in 
nearly all cases, the convergence is more variable.}
\end{figure}

\subsubsection{The diagonal case: targeting small eigenvalues}

\begin{example} \label{ex:10diag}
Consider the diagonal matrix
\[ A = {\rm diag}(0.01,\, 0.1,\, 3,\, 4,\, 5,\, \ldots,\, 9,\, 10).\]
Here $b=r_0$ is a Normal(0,1) vector (MATLAB's {\tt randn} 
with {\tt rng(1319)}), scaled so $\|r_0\|_2=1$.
\Cref{fig:gwf2a} compares W-GMRES($m$) to GMRES($m$) for $m=3$ and $m=6$:
in both cases, weighting gives a big improvement.  \Cref{fig:hist} compares
convergence of GMRES(3) and W-GMRES(3) for 
5000~Normal(0,1) right-hand side vectors (also scaled so $\|r_0\|_2=1$).  
As in \Cref{fig:gwf2a}, W-GMRES(3) usually substantially improves convergence, 
but it also adds considerable variation in performance.  
(We have noticed this tendency overall: weighting seems to add 
more variation, depending on the right-hand side and finite-precision effects.)

Focus on the $m=3$ in \cref{fig:gwf2a}.  GMRES(3) converges slowly: 
the small eigenvalues of $A$ make this problem difficult. 
\Cref{fig:gwf2bc} shows the residual components in each eigenvector 
after cycles~2, 3 and~4: GMRES(3) makes little progress in the smallest eigencomponent. 
Contrast this with W-GMRES(3):  after cycles~2 and~3 all eigencomponents 
have been reduced \emph{except the first one, corresponding to $\lambda=0.01$}.
Cycle~4 emphasizes this neglected eigenvalue,
reducing the corresponding eigencomponent in the residual by two orders of magnitude
(while increasing some of the other eigencomponents).  Without weighting, GMRES(3)
does not target this eigenvalue.

\begin{figure}
\vspace{.10in}
   \begin{center}
      \includegraphics[scale=0.48]{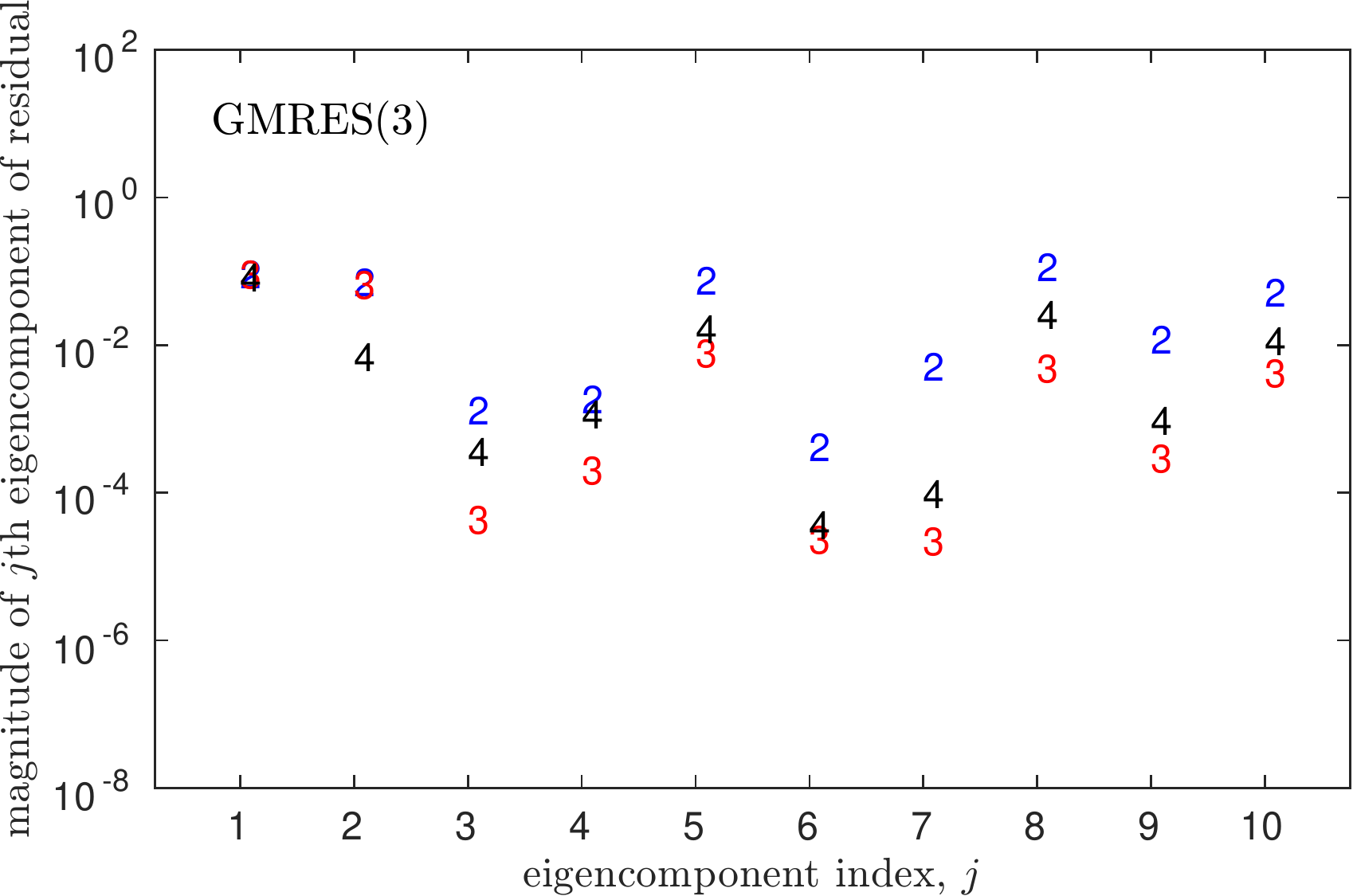} 
      \begin{picture}(0,0)
         \put(-195.7,50){\vector(0,1){15}}
         \put(-206,40){\footnotesize $\lambda=0.01$}
      \end{picture}

\vspace*{.5em}
      \includegraphics[scale=0.48]{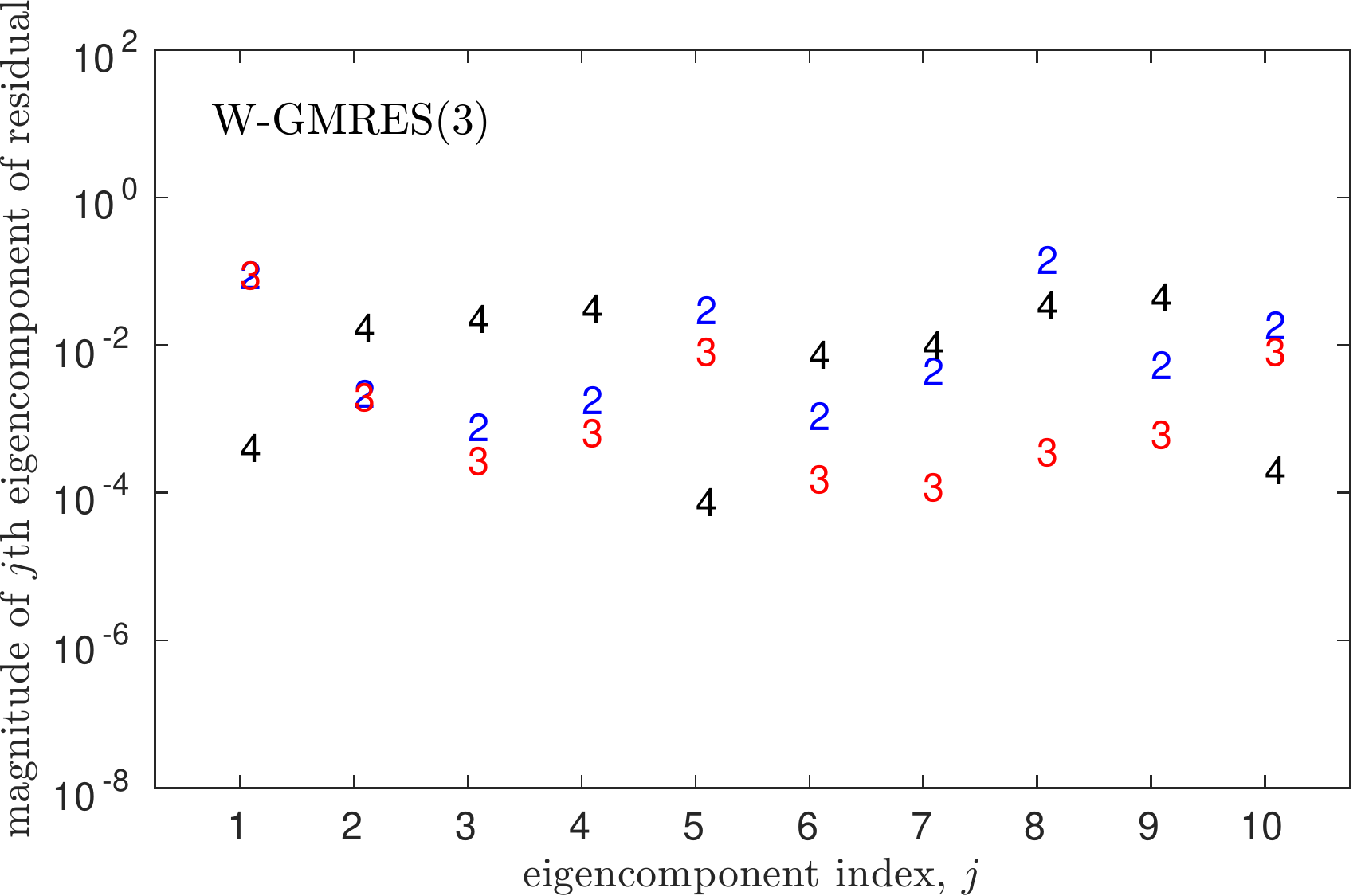} 
      \begin{picture}(0,0)
         \put(-195.7,50){\vector(0,1){15}}
         \put(-206,40){\footnotesize $\lambda=0.01$}
      \end{picture}
   \end{center}

\vspace*{-6pt}
\caption{\label{fig:gwf2bc}
Eigencomponents of the residual at the end of cycles
$\ell=2$, $3$, and $4$ for standard GMRES(3) (top) and W-GMRES(3) (bottom)
applied to the matrix in Example~2.  Cycle~4 of W-GMRES(3)
reduces the eigencomponent corresponding to $\lambda_1 = 0.01$ by two
orders of magnitude.}
\end{figure}
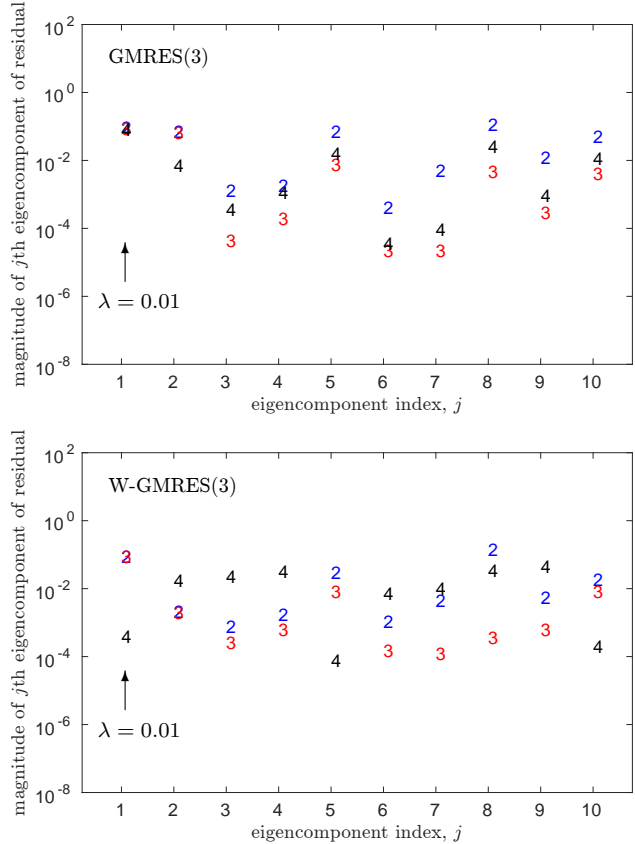

\Cref{fig:ex2poly} examines the GMRES residual polynomials $\poly_\ell$ 
that advance GMRES($m$) and W-GMRES($m$) from cycle~$\ell$ to cycle~$\ell+1$ 
(i.e., $r_{(\ell+1)m} = \poly_\ell(A) r_{\ell m}$ with $\poly_\ell(0)=1$)
for $\ell=2,\ldots, 5$.
GMRES(3) shows a repeating pattern,
with every other polynomial very similar;  
such behavior has been noted before~\cite{BaJeMa,CC}.  
None of these four polynomials is small at $\lambda=0.01$;
see the zoomed-in image on the right of \cref{fig:ex2poly}.  
(At cycle~4, $\poly_4(0.01) \approx 0.908$.)
The W-GMRES(3) polynomials at cycles~2, 3, and~5 are much like those for GMRES(3);
however, the polynomial for cycle~4 breaks the pattern, 
causing a strong reduction in the first eigencomponent:  
$\poly_4(0.01)\approx 0.00461$. 
Note that $\poly_4$ is also small at $\lambda=5$ and $\lambda=10$.
These are the other eigencomponents that are reduced at this cycle; the others
all increase, as seen in \cref{fig:gwf2bc}.
While the weighted norm decreases from~0.0906 to~0.0061 during cycle~4, 
the 2-norm shows a more modest reduction, from~0.0913 to~0.0734.
Indeed, experiments with different initial residuals often produce similar results, though the cycle that targets $\lambda=0.01$ can vary, and the 2-norm of the W-GMRES residual can increase at that cycle.

\begin{figure}
\vspace{.10in}
\begin{center}
   \includegraphics[scale=0.42]{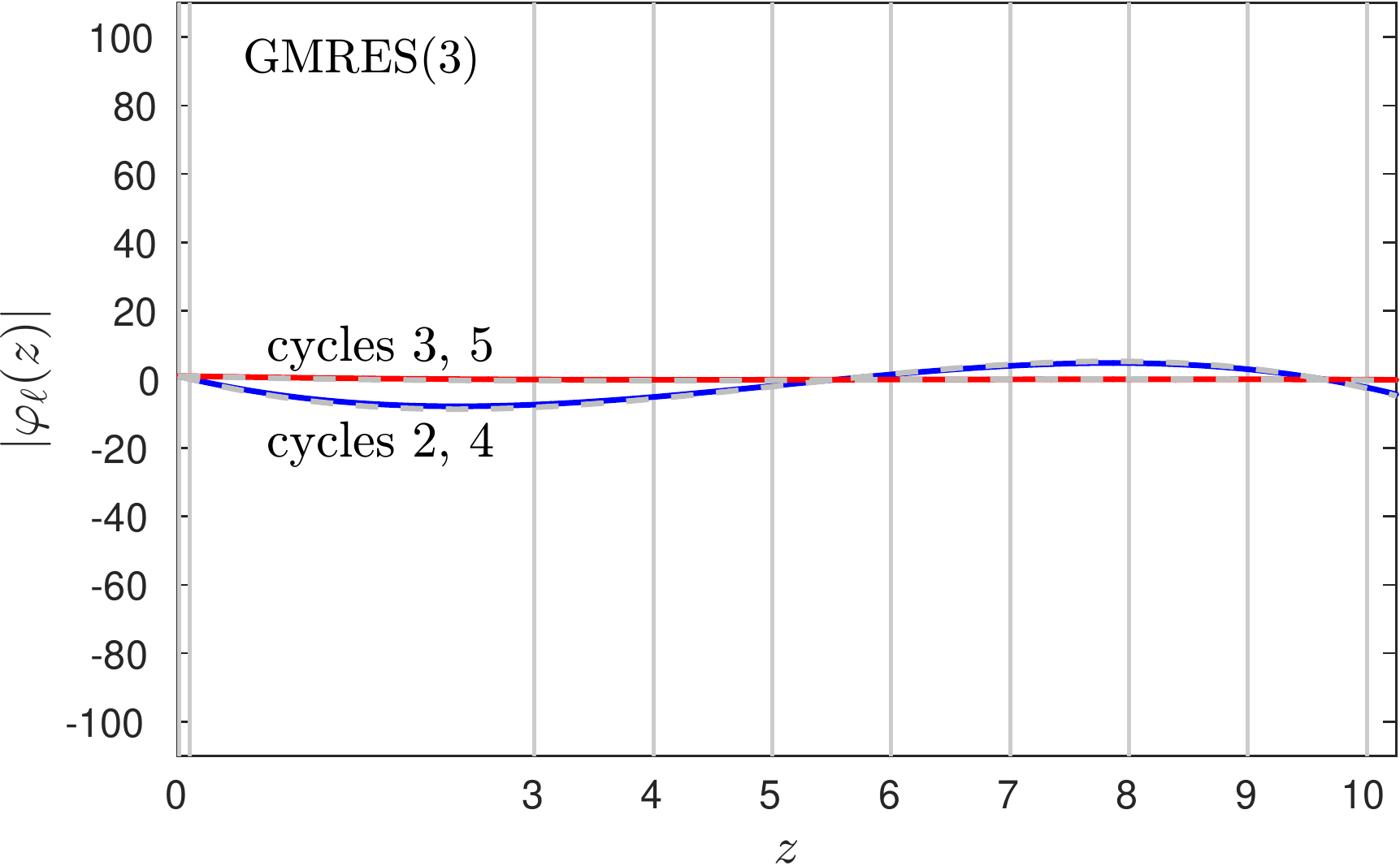}\quad
   \includegraphics[scale=0.42]{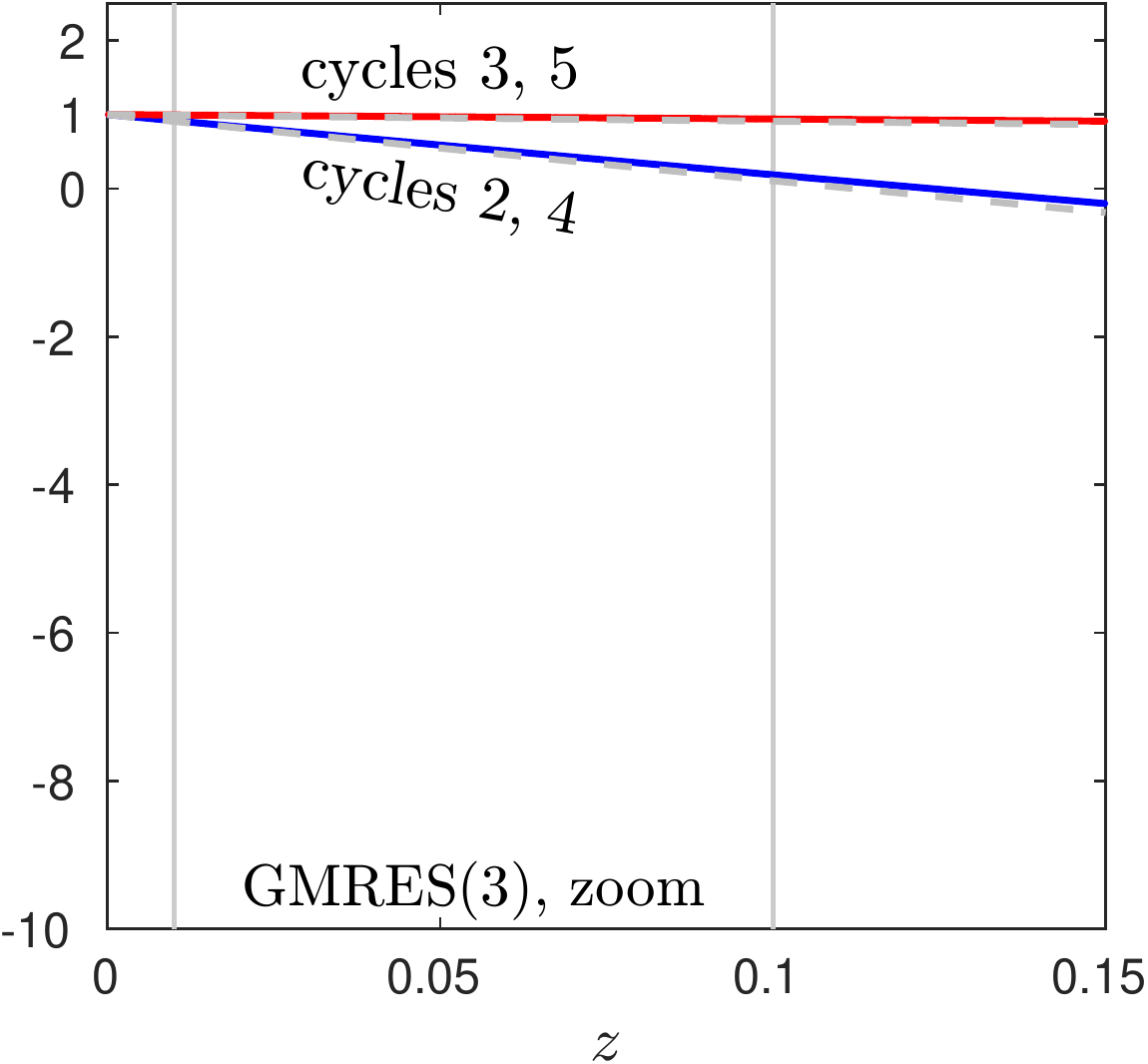} 

\vspace*{5pt}
   \includegraphics[scale=0.42]{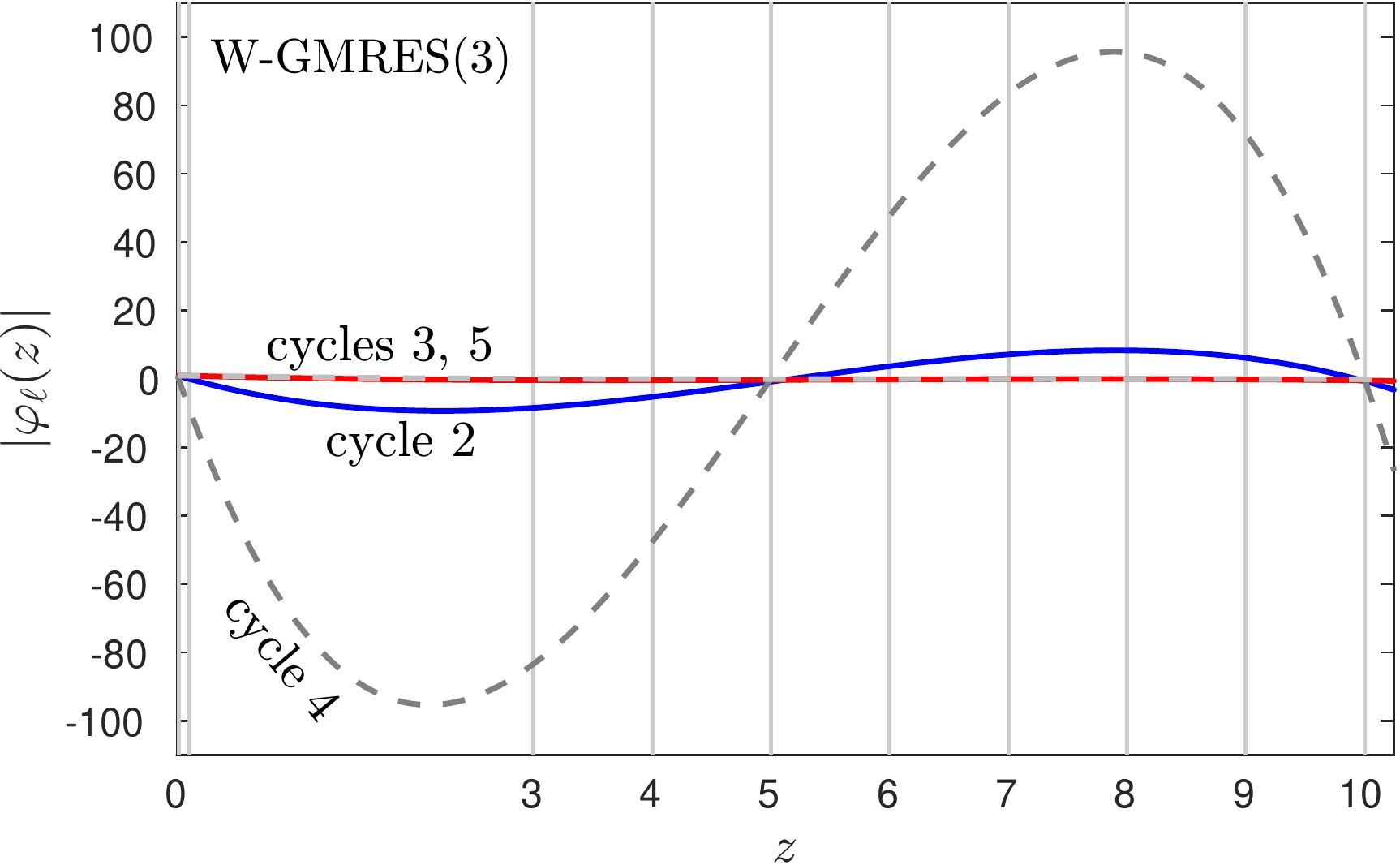}\quad
   \includegraphics[scale=0.42]{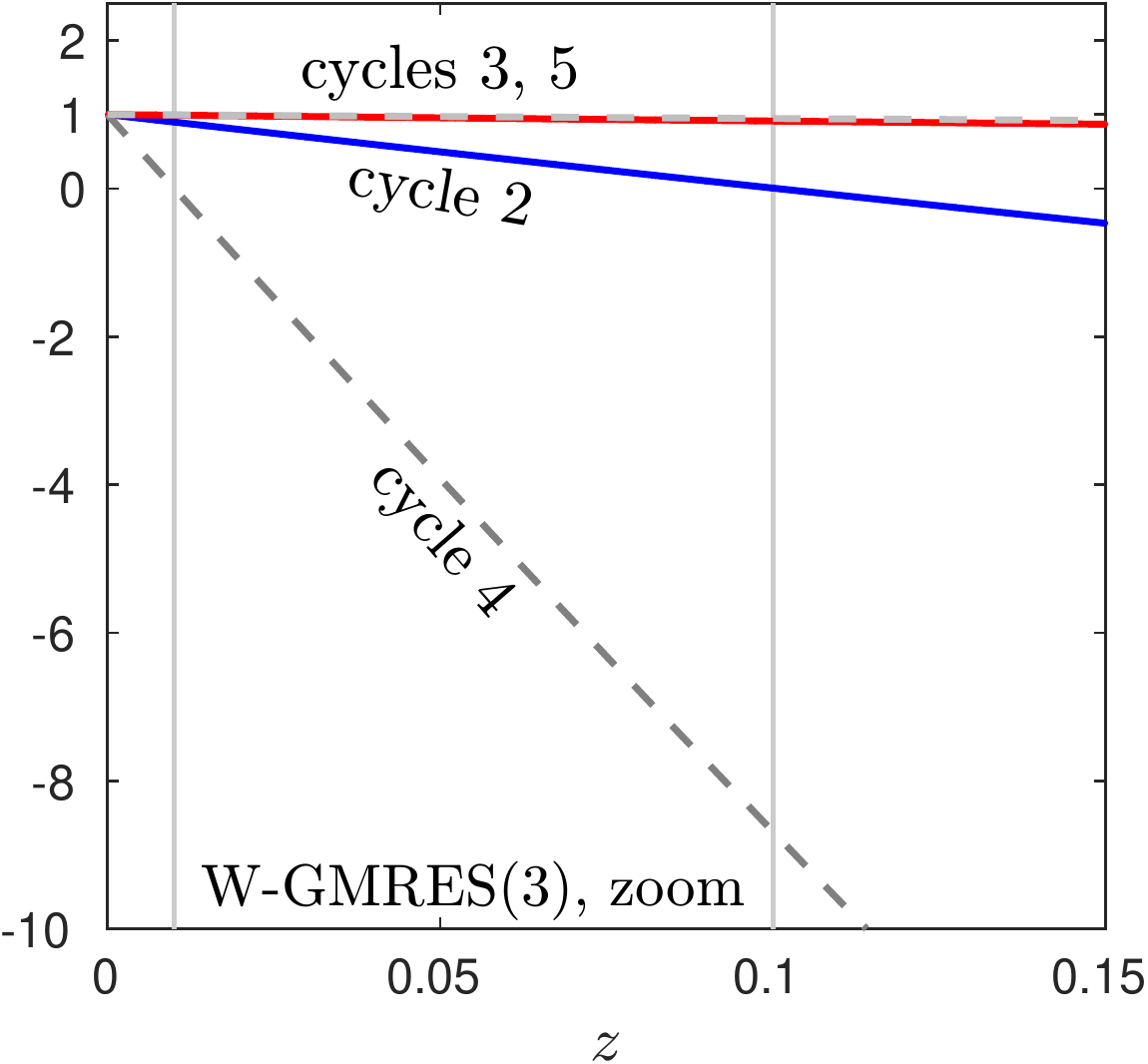} 
\end{center}

\vspace*{-8pt}
\caption{ \label{fig:ex2poly}
Residual polynomials for cycles 2--5 of standard GMRES(3) (top)
and W-GMRES(3) (bottom) for \cref{ex:10diag}.  The gray vertical lines
show the eigenvalues of $A$.  At cycle~4, the 
weighted inner product allows W-GMRES(3) to target the
first eigencomponent, increasing other components.
}
\vspace*{-4pt}
\end{figure}
\end{example}

The next result makes such observations precise for $2\times 2$ matrices.

\begin{theorem} \label{thm:g1roots}
Let $A = {\rm diag}(\lambda, \ 1)$ 
with real $\lambda\ne 0$, and suppose the residual vector at the 
start of a cycle of GMRES is
$b = [\thinspace b_1 \ \thinspace b_2 \thinspace]^T$
for $b_1, b_2 \in \R$ with $b_1 \ne 0$.
Define $\beta := b_2/b_1$.
Then the GMRES(1) polynomial has the root 
\[   \frac {\lambda^2 + \beta^2 } {\lambda + \beta^2},\]  
while the W-GMRES(1) polynomial has the root
\[   \frac {\lambda^2 + |\beta|^3} {\lambda + |\beta|^3}.\]  
\end{theorem}

The proof is a straightforward calculation.  
To appreciate this theorem, take $\lambda = \beta = 0.1$, giving the roots  
0.1818 for GMRES(1) and~0.1089 for W-GMRES(1).  
These polynomials cut the residual component corresponding to $\lambda=0.1$
by a factor~0.450 for GMRES(1) and~0.0818 for W-GMRES(1), the latter causing
W-GMRES(1) to increase the component for $\lambda=1$ by a factor of~8.18.  
This increase is worthwhile, given the reduction in the tough small component; 
later cycles will handle the easier component.
For another example, take $\lambda = \beta = 0.01$: 
GMRES(1) reduces the $\lambda=0.01$ component by~0.495, while W-GMRES(1) 
reduces it by~0.0098, two orders of magnitude.  

\subsubsection{The diagonal case: breaking patterns in residual polynomials}

Beyond the potential to target small eigenvalues, W-GMRES can also break
the cyclic pattern into which GMRES($m$) residual polynomials often lapse.
Baker, Jessup and Manteuffel~\cite[Thm.~2]{BaJeMa} prove that if $A$ is an $n\times n$ 
symmetric (or skew-symmetric) matrix, then the GMRES($n-1$) 
residual vectors exactly alternate in direction, i.e., the residual 
vector at the end of a cycle is an exact multiple of the residual two cycles before;
thus the GMRES residual polynomial at the end of each cycle repeats the
same polynomial found two cycles before.  
Baker et al.~observe the same qualitative behavior 
for more frequent restarts, and suggest that disrupting this pattern 
(by varying the restart parameter~\cite{BaJeKo} 
or augmenting the subspace~\cite{BaJeMa})
can improve convergence. (Longer cyclic patterns can emerge for 
nonnormal $A$~\cite{CC}.)  
Changing the inner product can have a similar effect, 
with the added advantage of targeting difficult eigenvalues if the corresponding
eigenvectors are well-disposed.
While Essai's residual-based weighting can break cyclic patterns, other schemes 
(such as random weighting~\cite{GP14,SNZa}; see \cref{sec:randw})
can have a similar effect; 
however, such arbitrary inner products take no advantage 
of eigenvector structure.
The next example cleanly illustrates how W-GMRES can break patterns.

\begin{figure}[b!]
\begin{center}
\includegraphics[scale=0.35]{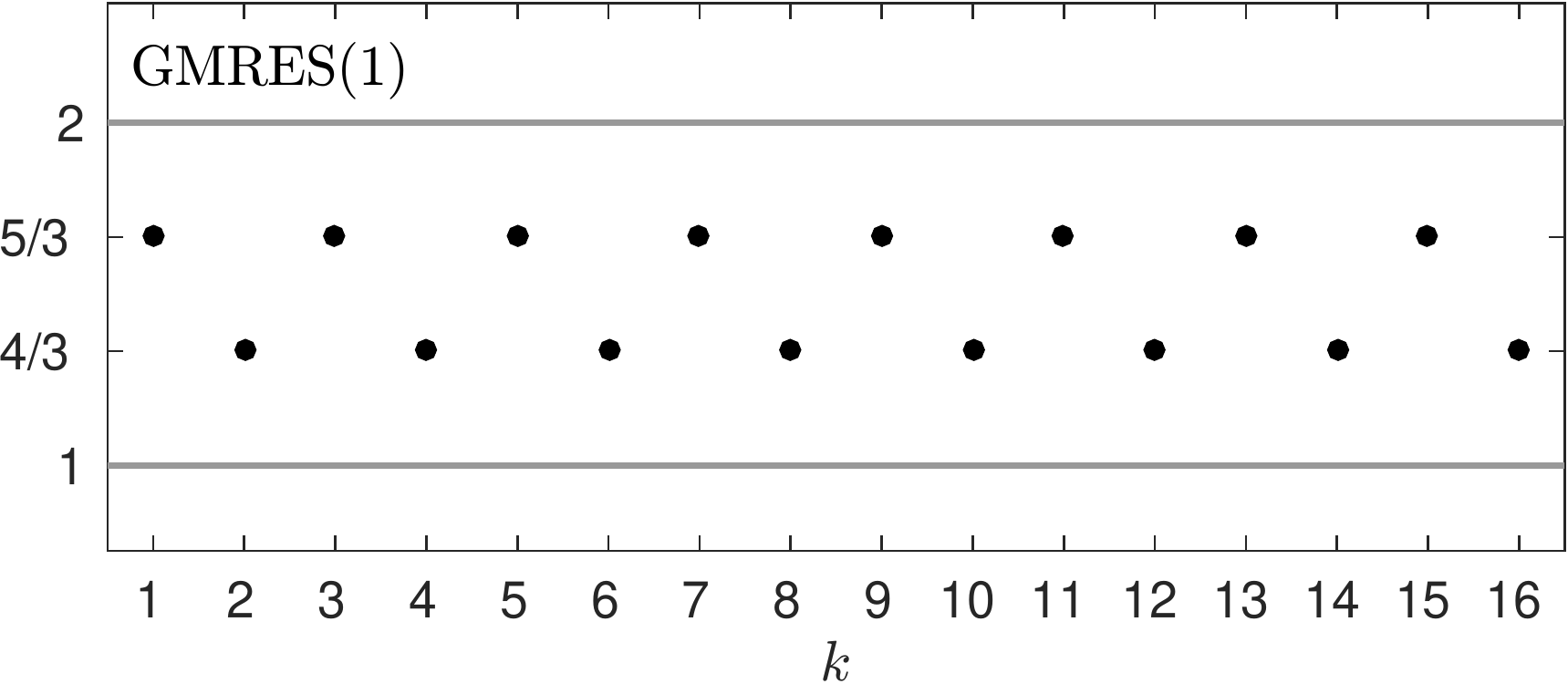}
\quad
\includegraphics[scale=0.35]{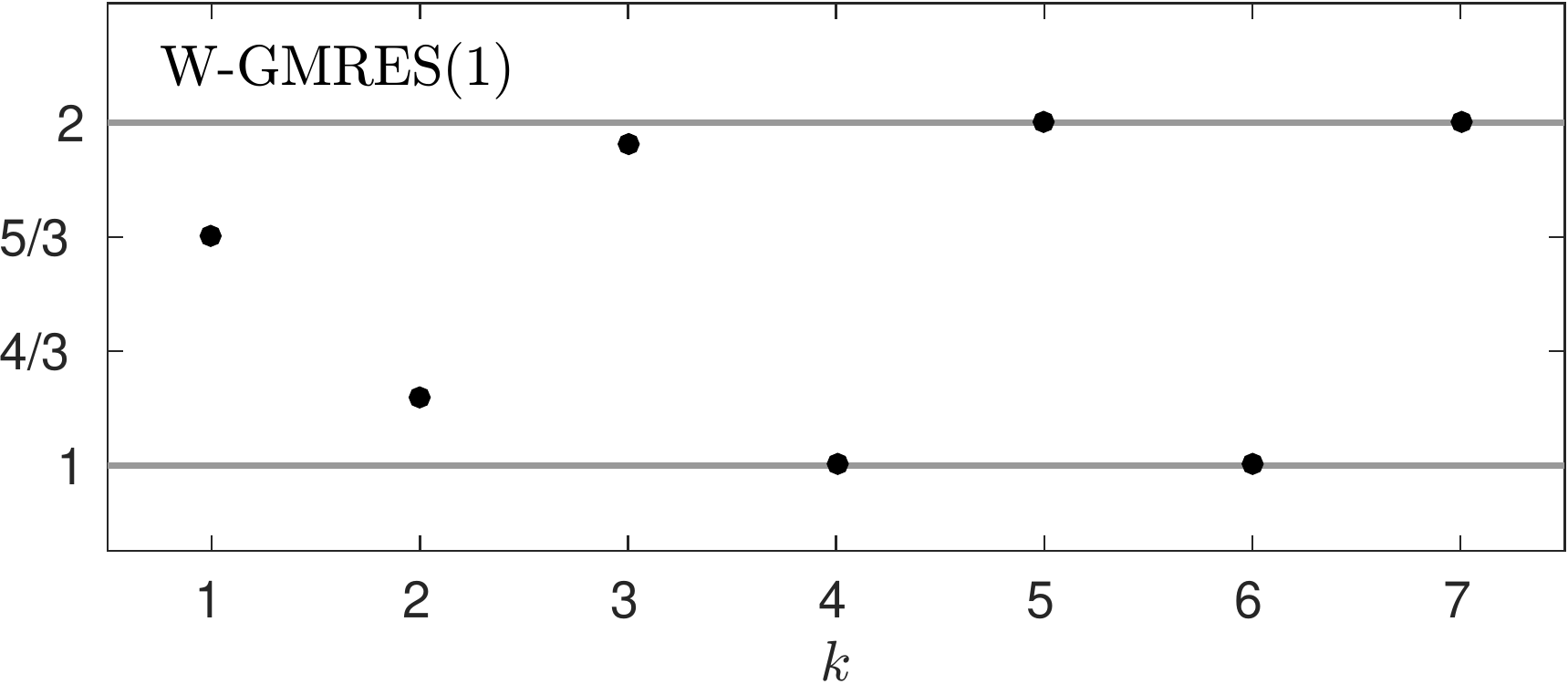}
\end{center}

\vspace*{-8pt}
\caption{\label{fig:g1roots}
Roots of the GMRES(1) and W-GMRES(1) residual polynomials for $A = {\rm diag}(2,1)$ 
with $b=[1, 1]^T$, as a function of the restarted GMRES cycle, $k$.  
The GMRES(1) roots occur in a repeating pair, and the method
takes 16~iterations to converge to $\|r_k\|_2/\|r_0\|_2 \le 10^{-8}$.  In contrast, 
the W-GMRES(1) roots are quickly attracted to the eigenvalues, giving 
convergence in just 7~steps.}
\end{figure}

\begin{example} \label{ex3}
Apply GMRES(1) to $A = {\rm diag}(2,1)$ with $b = [\thinspace 1 \ \thinspace 1 \thinspace]^T$.  
The roots of the (linear) GMRES(1) residual polynomials alternate between $5/3$ and $4/3$ (exactly), 
as given by~\cref{thm:g1roots}.
Therefore the GMRES(1) residual polynomials for the cycles alternate between 
$\poly_k(z) = 1 - (3/5) z$ and $\poly_k(z) = 1 - (3/4) z$.  
GMRES(1) converges ($\|r_k\|_2/\|r_0\|_2 \le 10^{-8}$) in 16~iterations.
W-GMRES(1) takes only~7 iterations.  
As \cref{fig:g1roots} shows, weighting breaks the cycle of polynomial roots, 
giving 1.667, 1.200, 1.941, 1.0039, 1.999985, 1.0000000002, and $\approx 2$:  
the roots move out toward the eigenvalues, alternating which eigenvalue 
they favor and reducing $\|r_k\|_2$ much faster.
\end{example}

\subsection{Measuring eigenvector localization}  \label{sec:local}

The justification for W-GMRES in~(\ref{eq:diag3}) required diagonal $A$,
so the eigenvectors are columns of the identity matrix.  
When $A$ is not diagonal, the motivation for W-GMRES is less compelling.
Yet in many applications $A$ is not close to diagonal, but still
has \emph{localized} eigenvectors: only a few entries 
are large, and the eigenvectors resemble columns of the identity matrix.
Computational evidence suggests that W-GMRES
can still be effective, particularly when the eigenvectors associated 
with the small magnitude eigenvalues are localized.  
The Add20 matrix in \cref{ex:add20} is such an example.

We gauge the localization of the $p$ smallest magnitude 
eigenvectors by measuring how much they are concentrated
in their $p$ largest magnitude entries.%
\footnote{Certain applications motivate more specialized
measures of localization based on the rate of exponential
decay of eigenvector entries about some central entry.}
For diagonalizable $A$, label the eigenvalues in increasing magnitude, 
$|\lambda_1|\le |\lambda_2| \le \cdots \le |\lambda_n|$,
with associated unit 2-norm eigenvectors $v_1, \ldots, v_n$.
For any $y\in \C^n$, let $y^{(p)} \in \C^p$ denote the subvector
containing the $p$ largest magnitude entries of $y$ (in any order).
Then 
\begin{equation} \label{eq:locp}
 \loc_p(A) := {1\over \sqrt{p}} \Big(\sum_{j=1}^p {\|v_j^{(p)}\|_2^2}\Big)^{1/2}.
\end{equation}

\begin{figure}[t!]
\begin{picture}(0,155)
 \put(0,0){\includegraphics[scale=0.5]{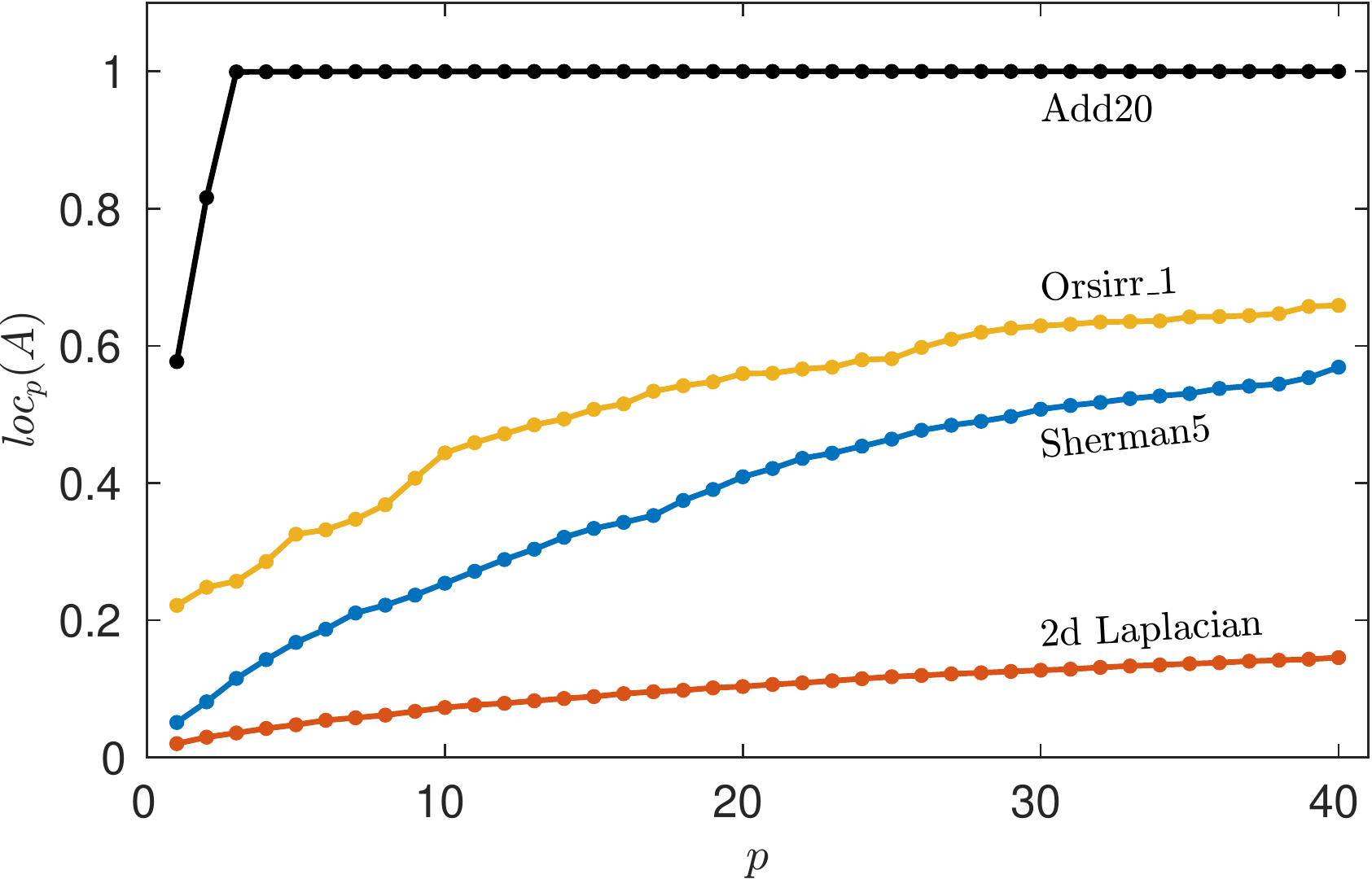}}
 \put(245,42){\footnotesize \begin{tabular}{lcc}
             \emph{matrix} &  $n$ & $\loc_5(A)$ \\ \hline
              Add20 & $2395$ & $0.9996$ \\ 
              Orsirr\_1& $1030$ & $0.3253$\\ 
              Sherman5 & $3312$ & $0.1676$\\ 
              2d Laplacian\!\!\! & $9801$ & $0.0475$ \\ 
             \end{tabular}
            }
\end{picture}

\vspace*{-7pt}
\caption{\label{fig:comploc}
The localization measures $\loc_p(A)$ for $p=1,\ldots, 40$ for four matrices
that are used in our W-GMRES experiments.
If $\loc_p(A)$ is close to~1, the eigenvectors associated with smallest 
magnitude eigenvalues are highly localized.}
\end{figure}

\noindent
Notice that $\loc_p(A) \in(0,1]$, and $\loc_n(A) = 1$.  
If $\loc_p(A) \approx 1$ for some $p\ll n$, the eigenvectors associated with
the $p$ smallest magnitude eigenvalues are localized within $p$ positions.
This measure is imperfect, for it does not include the magnitude of the
eigenvalues nor potential non-orthogonality of the eigenvectors.
Still $\loc_p(A)$ can be a helpful instrument for assessing localization.  
\Cref{fig:comploc} shows $\loc_p(A)$
for $p=1,\ldots, 40$ for four matrices we use in our experiments.%
\footnote{Two of these matrices, Sherman5 and the 2d Laplacian, have repeated
(derogatory) eigenvalues; according to the principle of \emph{reachable
invariant subspaces} (see, e.g., \cite{BeEmRo,LS13}), the GMRES process 
only acts upon one eigenvector in the invariant subspace,
given by the component of the initial residual in that space.  
Thus the localization measures in \cref{fig:comploc} 
include only one eigenvector for the derogatory eigenvalues,
specified by the initial residual vectors used in our experiments.}
The Add20 matrix has strongly localized eigenvectors; 
the discrete 2d~Laplacian has eigenvectors that are far from
localized, i.e., they are \emph{global}.  The \Orsirr\ and Sherman5
matrices have eigenvectors with intermediate localization.

\subsection{The nondiagonal case}

When $A$ is not diagonal, the justification for W-GMRES in 
equation~(\ref{eq:diag3}) is lost; the effect of weighting becomes more subtle.
The weights can still break cyclic patterns in restarted GMRES, but the interplay
between the weights and the eigenvectors is more difficult to understand.
Suppose $A$ is diagonalizable, $A = V \Lambda V^{-1}$, still with 
$W = S^*S$.  
Unlike the diagonal case, $S$ will not generally commute with $V$ and $\Lambda$, 
so~\cref{eq:diag} is replaced by 
\begin{subequations} \label{eq:nondiag}
\begin{eqnarray} 
 \|r_k\|_\W  &=& \mingmres \|\poly(SAS^{-1})Sb\|_2 
                   \label{eq:nondiag1} \\
               &=& \mingmres \|(SV) \poly(\Lambda )(SV)^{-1} Sb\|_2. 
                   \label{eq:nondiag2} \\[-1em] \nonumber
\end{eqnarray}
\end{subequations}
The matrix $S$ \emph{transforms the eigenvectors of $A$}~\cite{GL01}.  
Suppose, impractically, that $V^{-1}$ were known, allowing the \emph{nondiagonal} weight
\[S = {\rm diag}(s_1,\ldots, s_n) V^{-1}.\]
In this case~(\ref{eq:nondiag2}) reduces to 
\begin{equation} \label{eq:evweight}
 \|r_k\|_\W^2 = \mingmres \sum_{j=1}^n  |@\poly(\lambda_j)|^2 |s_j b_j|^2,
\rlap{\kern18pt\emph{eigenvector weighting}}
\end{equation}
a perfect analogue of the diagonal case~(\ref{eq:diag3}) that could appropriately
target the small magnitude eigenvalues that delay convergence.
One might \emph{approximate} this eigenvector weighting 
by using, as a proxy for $V^{-1}$, some estimate of the left eigenvectors of $A$ associated
with the smallest magnitude eigenvectors.  (The rows of $V^{-1}$ are
the left eigenvectors of $A$, since $V^{-1} A = \Lambda V^{-1}$.)
We do not pursue this idea here, instead using eigenvector information 
via deflated restarting in \cref{sec:drwgmres}.

The conventional W-GMRES($m$) algorithm instead uses diagonal $W$ (and $S$).
In this case, $SV$ scales the \emph{rows} of $V$, effectively emphasizing
certain entries of the eigenvectors at the expense of others.  
Write out the right and left eigenvectors,
\[ V = [v_1\ v_2\ \cdots\ v_n] \in \C^{n\times n}, \qquad
 V^{-1} = \left[\begin{array}{c} \wh{v}_1^* \\ \wh{v}_2^* \\ \vdots \\ \wh{v}_n^* \end{array}\right]\in\C^{n\times n}.\]
Let $c := V^{-1} b$ denote the expansion coefficients for $b$ in the 
eigenvectors of $A$: $b= Vc$.
Then one can also render~(\ref{eq:nondiag2}) in the form
\begin{equation}
 \|r_k\|_\W  = \mingmres \Big\| \sum_{j=1}^n  c_j @@ \poly(\lambda_j) @ Sv_j \Big\|_2,
\rlap{\kern27pt \emph{diagonal weighting}}
\end{equation}
a different analogue of~(\ref{eq:diag3}).
Denote the $\ell$th entry of $v_j$ by $(v_j)_\ell$.
In most cases $V$ will be dense.
Suppose that only the $q$th entry of $r_0$ is large, and that $(v_j)_q \ne 0$
for $j=1,\ldots, n$.
Then Essai's weighting makes $|(Sv_j)_q| \gg |(Sv_j)_\ell|$ for all $\ell \ne q$:
all the vectors $Sv_j$ form a small angle with the $q$th column of the identity.
Thus $SV$ is ill-conditioned and $SAS^{-1}$ will have a large departure from normality, 
often a troublesome case for GMRES; see, e.g.,~\cite[chap.~26]{TE05}.
It is not evident how such weighting could help W-GMRES($m$) 
focus on small eigenvalues, as it does for diagonal~$A$,
suggesting an explanation for the mixed performance of W-GMRES($m$)~\cite{CaYu,GP14}.

\subsubsection{\boldmath Stagnation of W-GMRES($m$) when GMRES($m$) converges}
We show that by transforming the eigenvectors, diagonal weighting can even prevent
W-GMRES($m$) from converging at all.  Faber et al.~\cite{FJKM96} prove that GMRES($m$) converges
for all initial residuals provided there exists no $v\in\C^n$ such that
\[ v^*A^k v = 0, \quad \mbox{for all $k=1,\ldots, m$.}\]
Thus GMRES($m$) will converge for all $m\ge 1$ and $r_0$ if the field of values of $A$,
\[ F(A) := \{ v^*Av : v\in\C^n, \|v\|_2 = 1\},\]
does not contain the origin.  
Weighting effectively applies GMRES to the transformed matrix $SAS^{-1}$,
and now it is possible that $0 \in F(SAS^{-1})$ even though $0 \not \in F(A)$.
In extreme cases this means that W-GMRES($m$) can stagnate even when GMRES($m$) converges,
as shown by the next example.
(In contrast, for $A$ and $S$ diagonal, $SAS^{-1} = A$, so $F(A)=F(SAS^{-1})$.)

\begin{figure}[b!]
\begin{center}\includegraphics[scale=0.4]{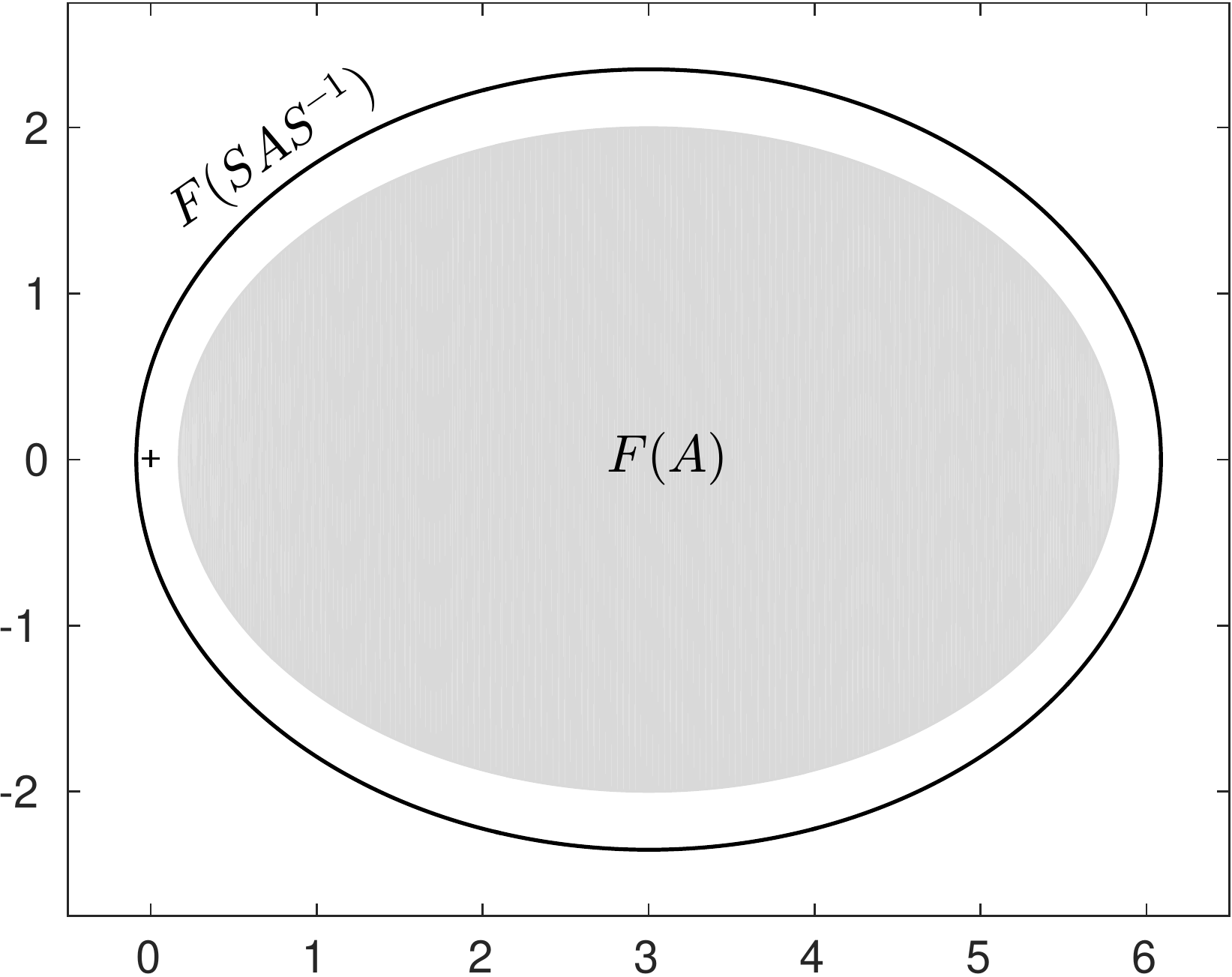}\end{center}

\vspace*{-3pt}
\caption{\label{fig:stag2fv}
The fields of values $F(A)$ (gray region) and $F(SAS^{-1})$ (boundary 
is a black line) in $\C$
for $A$ and $r_0$ in~\cref{eq:stag2}; the cross ($+$) marks the origin.
Since $0\not\in F(A)$, GMRES($m$) converges for any $r_0$;
however, $0\in F(SAS^{-1})$ and $(Sr_0)^*(SAS^{-1})(Sr_0)=0$, so W-GMRES(1)
completely stagnates.}
\end{figure}

\begin{example}[W-GMRES(1) stagnates while GMRES($1$) converges] \label{ex:stagnate}\ \\
Consider 
\begin{equation} \label{eq:stag2}
 A = \left[\begin{array}{rr} 1 & -4 \\ 0 & 5 \end{array}\right], \qquad
   r_0 = \left[\begin{array}{c} 1 \\ {1\over 10}(5+\sqrt{5})\end{array}\right]
     = \left[\begin{array}{c} 1 \\ 0.72360\ldots\end{array}\right].
\end{equation}
\Cref{fig:stag2fv} shows $F(A)$, an ellipse in the complex plane.  
From the extreme eigenvalues of the Hermitian and skew-Hermitian parts of $A$,
one can bound $F(A)$ within a rectangle in $\C$ that does not contain the origin,
\[ {\rm Re}(F(A)) = [3-2\sqrt{2}, 3+2 \sqrt{2}] \approx [0.17157,5.82843],
   \qquad
   {\rm Im}(F(A)) = [-2@{\rm i}, 2@{\rm i}@],\]
and hence GMRES(1) converges for all initial residuals.
Now consider the specific $r_0$ given in~\cref{eq:stag2}.
Using the scaling $S = {\rm diag}(\sqrt{|(r_0)_1|/\|r_0\|_\infty}, \sqrt{|(r_0)_2|/\|r_0\|_\infty})$
at the first step gives the W-GMRES(1) polynomial (see, e.g., \cite[eq.~(2.2)]{Em03})
\[  r_1 = \poly(A)r_0 = r_0 - {(Sr_0)^*(SAS^{-1})(Sr_0) \over \| (SA S^{-1}) (Sr_0)\|_2^2 } (SAS^{-1}) (Sr_0) 
              = r_0 - {r_0^* S^2 A r_0 \over \|SAr_0\|_2^2} S A r_0. \]
This example has been engineered so that $0\in F(SAS^{-1})$; indeed, 
\[ {\rm Re}(F(SAS^{-1})) = \bigg[3 - \sqrt{14-2\sqrt{5}}, 
               3+\sqrt{14-2\sqrt{5}}\,\bigg] \approx [-0.086724,6.08672],\]
as seen in \cref{fig:stag2fv}.
Worse still, this $r_0$ gives $r_0^* S^2 A r_0 = 0$, so $r_1=r_0$:
W-GMRES(1) makes no progress, and the same weight is chosen for the next cycle.
The weighted algorithm completely stagnates, even though GMRES(1) 
converges for any $r_0$.
\end{example}

\begin{figure}[b!]
\begin{center} \includegraphics[scale=0.5]{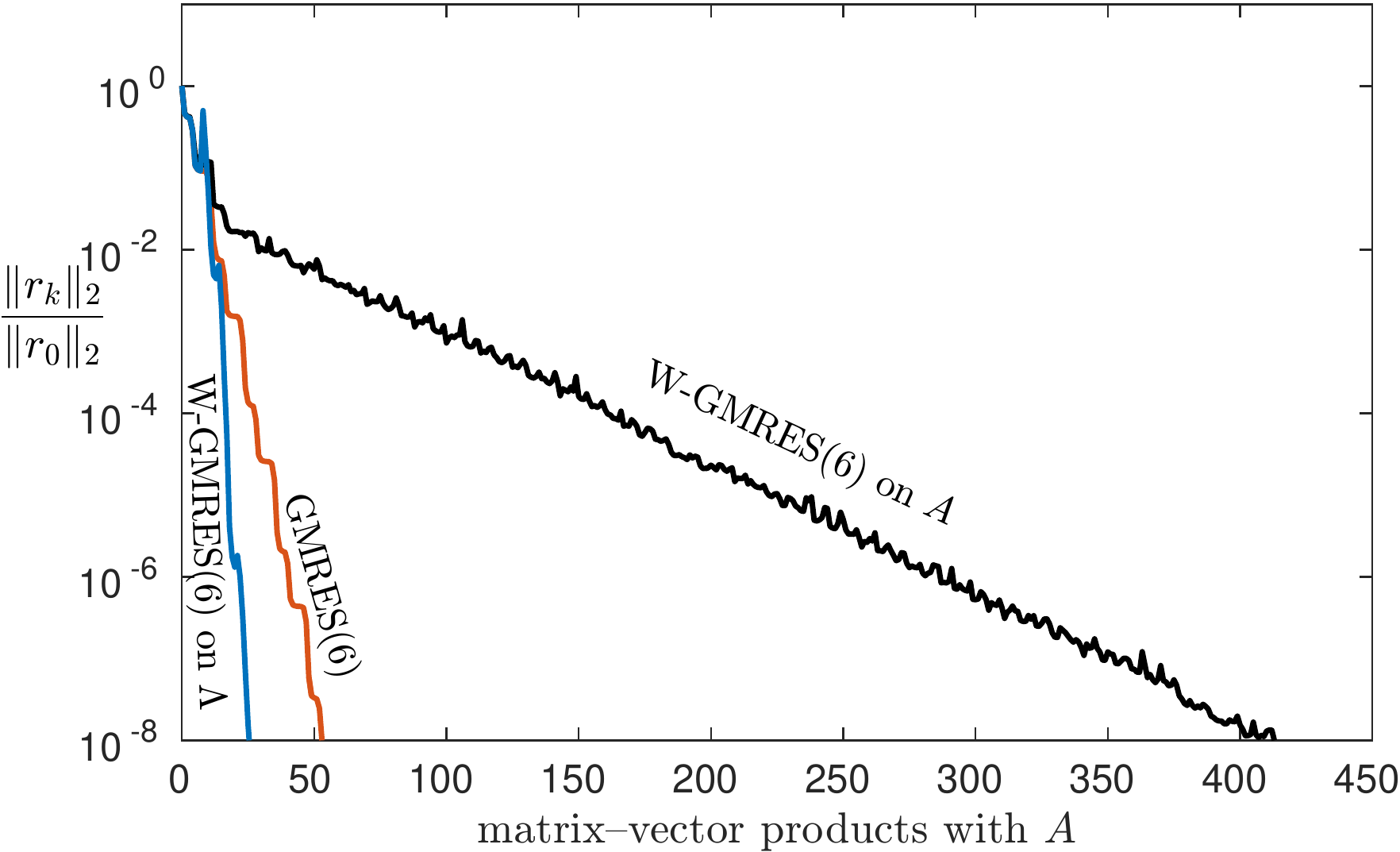} \end{center}

\vspace*{-3pt}
\caption{\label{fig:lapev}
GMRES(6) and W-GMRES(6) for \cref{ex:10lap}.
Applying GMRES(6) to $(\Lambda, b)$ and $(A, Vb)$ produces identical residual 2-norms.
When applied to $(\Lambda, b)$, W-GMRES(6) converges faster than GMRES(6) 
(seen previously in \cref{fig:gwf2a});
when applied to $(A,Vb)$ the convergence is much slower
($A$ having non-localized eigenvectors).} 
\end{figure}

\subsubsection{\boldmath Poor performance when small eigenvalues have global eigenvectors}
The last example shows one way GMRES($m$) can outperform
W-GMRES($m$) when $A$ is not diagonal.
Next we show that global eigenvectors associated with small magnitude
eigenvalues can also contribute to poor W-GMRES($m$) convergence.
\medskip

\begin{example}[W-GMRES($m$) worse than GMRES($m$): global eigenvectors]\ \\ 
\label{ex:10lap} 
We revisit \cref{ex:10diag}, again with dimension $n=10$ and eigenvalues
\[ \Lambda = {\rm diag}(0.01, 0.1, 3, 4, 5, \ldots, 9, 10),\]
but now let the eigenvectors of $A$ equal those of the 
symmetric tridiagonal matrix ${\rm tridiag}(-1,0,-1)$,
ordering the eigenvalues from smallest to largest.
(The $j$th column of $V$ has entries $(v_j)_\ell = \sin(j\ell \pi /(n+1))/\|v_j\|_2$;
see, e.g., \cite{Smi85}.)
The resulting $A=V\Lambda V^*$ is symmetric,
unitarily similar to $\Lambda$.
Let $b$ be the same Normal(0,1) vector used in \cref{ex:10diag}.
Since $V$ is unitary, GMRES($m$) applied to $(\Lambda, b)$ produces the same 
residual norms as GMRES($m$) applied to $(V\Lambda V^*, Vb)$.
However, weighted GMRES can behave very differently for the two problems:
when applied to $\Lambda$ (as in \cref{ex:10diag}), W-GMRES(6)
vastly outperforms GMRES(6).
For $A$ with its non-localized eigenvectors, W-GMRES(6) 
converges much more slowly, as seen in \cref{fig:lapev}.
\end{example}

The next example shows similar behavior for a larger matrix: 
eigenvalues on both sides of the origin and global eigenvectors give
poor W-GMRES($m$) performance.

\begin{figure}
  \begin{center}\includegraphics[scale=0.475]{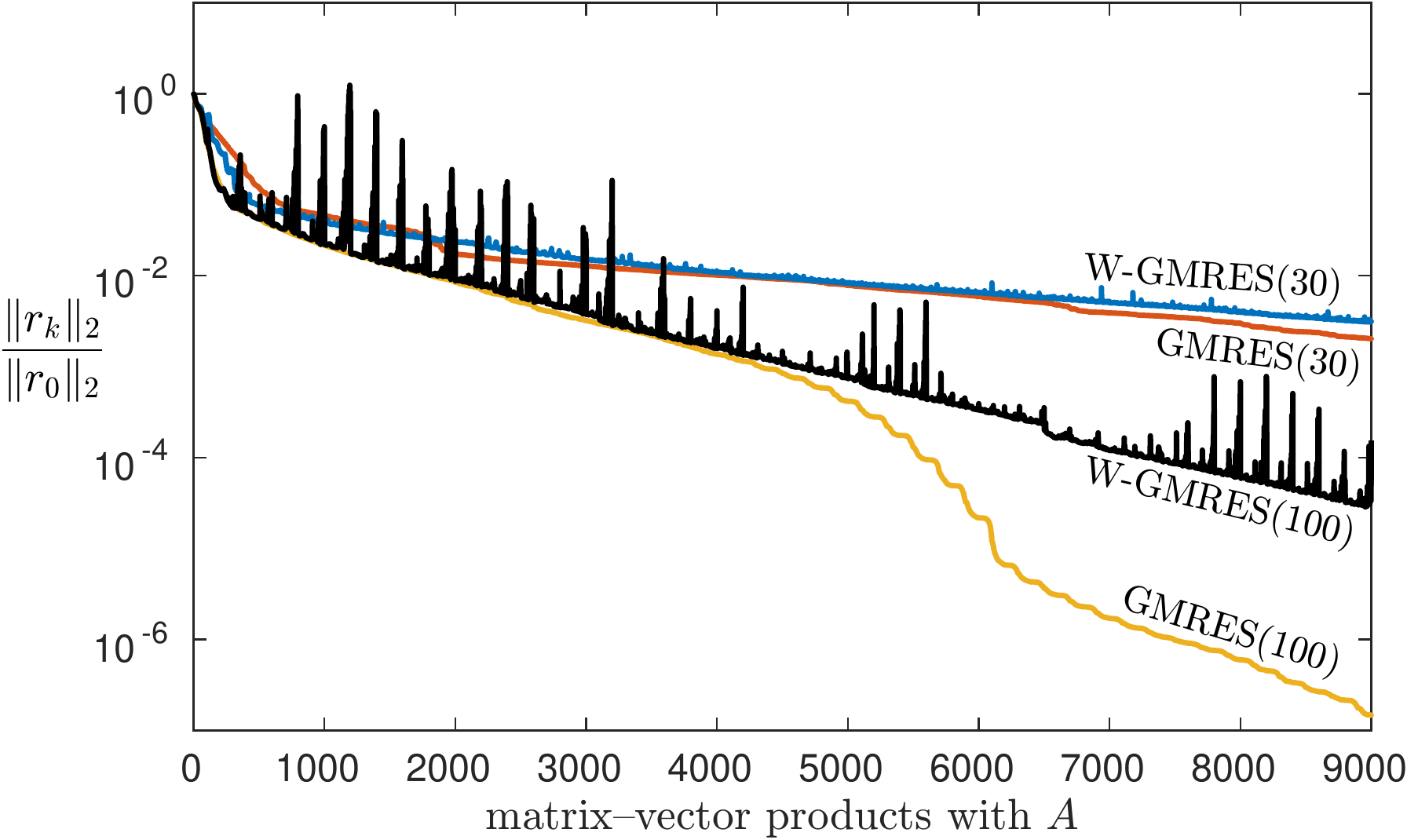} \end{center}

\vspace{4pt}
  \begin{center}\includegraphics[scale=0.475]{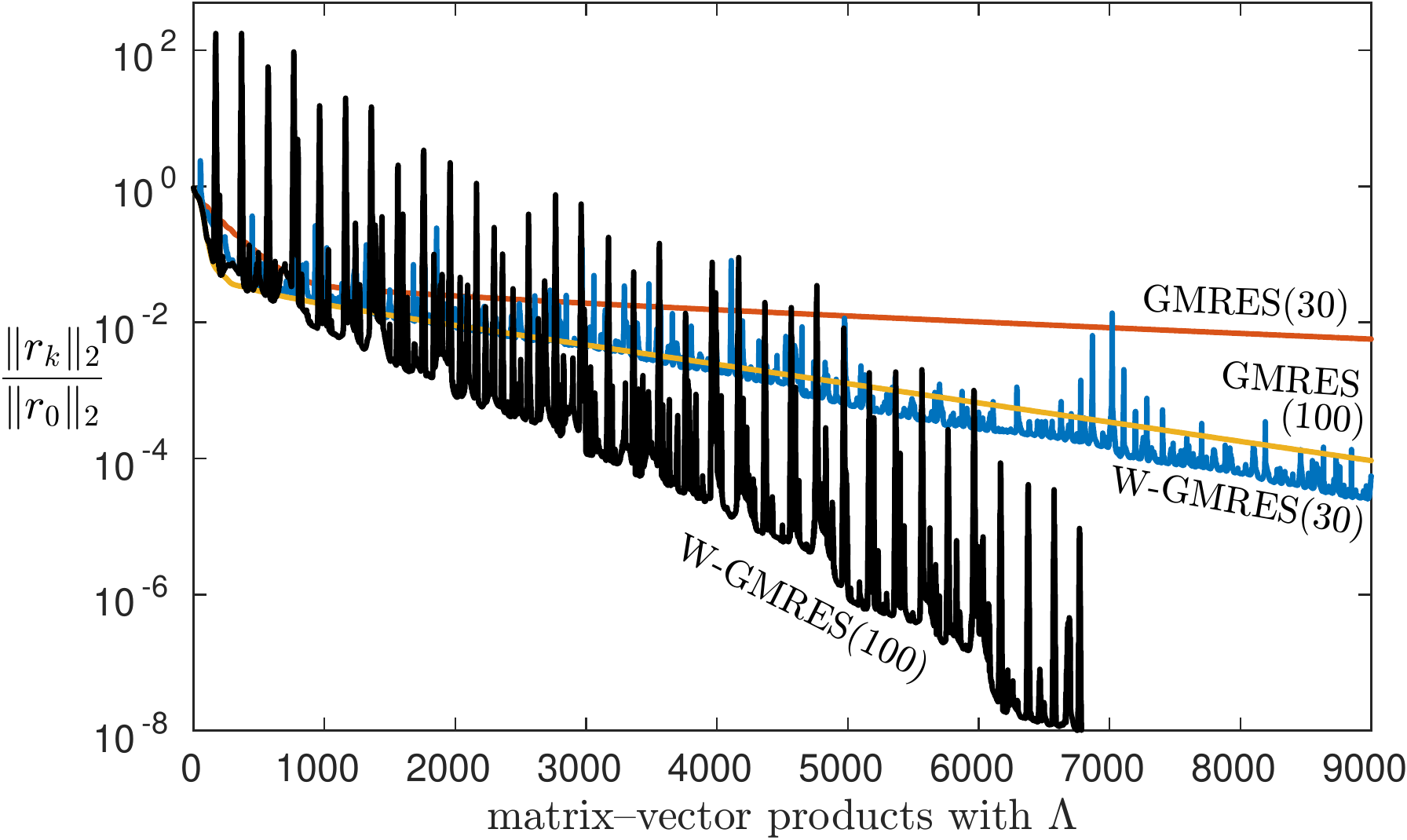} \end{center}

\vspace*{-4pt}
\caption{\label{fig:sherex}
GMRES($m$) and W-GMRES($m$) for the Sherman5 matrix.
On the top, the methods are applied to $Ax=b$, and 
weighting slows convergence.  
On the bottom, the same methods are applied to $\Lambda x = b$,
where $\Lambda$ is diagonal with the same eigenvalues of $A$; 
W-GMRES converges better.
(The GMRES(100) curve is yellow.)
For both problems, taking~$m=100$ steps 
with each weighted inner product gives large growth in $\|r_k\|_2$ 
between restarts.
}
\vspace*{-10pt}
\end{figure}

\begin{example}[W-GMRES($m$) slower than GMRES($m$): global eigenvectors]
\label{ex:sherman}\ \\ 
Consider the
nonsymmetric Sherman5 matrix from Matrix Market~\cite{MM}, with $n=3312$
and with $r_0$ a Normal(0,1) vector.
\Cref{fig:sherex}  shows that weighting slows restarted GMRES, 
particularly for the restart $m=100$.
\Cref{fig:comploc} shows that the eigenvectors are poorly localized,
e.g., in the unit eigenvector $v_1$ corresponding to the smallest eigenvalue,
the 100~largest components range only from 0.0509 to 0.0415 in magnitude.
This matrix has eigenvalues to the left and right of the origin, 
and the eigenvectors are not orthogonal.%
\footnote{The eigenvalue $\lambda=1$ has multiplicity 1674,
and the corresponding eigenvectors are columns of the identity.
The eigenvalues closest to the origin are well-conditioned, while some 
eigenvalues far from the origin are relatively ill-conditioned.  
(For the diagonalization computed by MATLAB's {\tt eig} command, 
$\|V\|_2 \|V^{-1}\|_2 \approx 3.237\times 10^6$.)
The pseudospectra of this matrix are well-behaved near the origin;
the effect of the eigenvalue conditioning on GMRES convergence is mild~\cite[chap.~26]{TE05}.}
The second plot in \cref{fig:sherex} shows how W-GMRES performance 
improves when $A$ is replaced by a diagonal matrix with the same spectrum.
\end{example}

\subsubsection{\boldmath Weighting sometimes helps despite global eigenvectors}
The past three examples suggest that W-GMRES performs poorly for 
matrices whose eigenvectors are not localized, but the algorithm is
more nuanced than that.  The next example is characteristic
of behavior we have often observed in our experiments.
We attribute the improvement over standard restarted GMRES to the
tendency of weighting to break cyclic patterns in the GMRES residual polynomials.

\begin{example}[W-GMRES($m$) better than GMRES($m$): global eigenvectors]
 \label{ex:lap2d} \ \\
Consider the standard 5-point finite difference discretization of the Laplacian on the unit square 
in two dimensions with Dirichlet boundary conditions. The uniform grid spacing $h=1/100$ 
in both directions gives a matrix of order $n=99^2 = 9801$; 
the initial residual is a random Normal(0,1) vector.
The small $\loc_p(A)$ values shown in \cref{fig:comploc} confirm 
that the eigenvectors are not localized;
in fact, they are less localized than for the Sherman5 matrix in the last example.
Despite this, \cref{fig:lap2d} shows that W-GMRES($m$) outperforms GMRES($m$)
for $m=10$ and $m=20$.
In this case, all the eigenvalues are positive, and none is particularly close to the origin.

\begin{figure}
\vspace{.10in}
\begin{center}\includegraphics[scale=0.5]{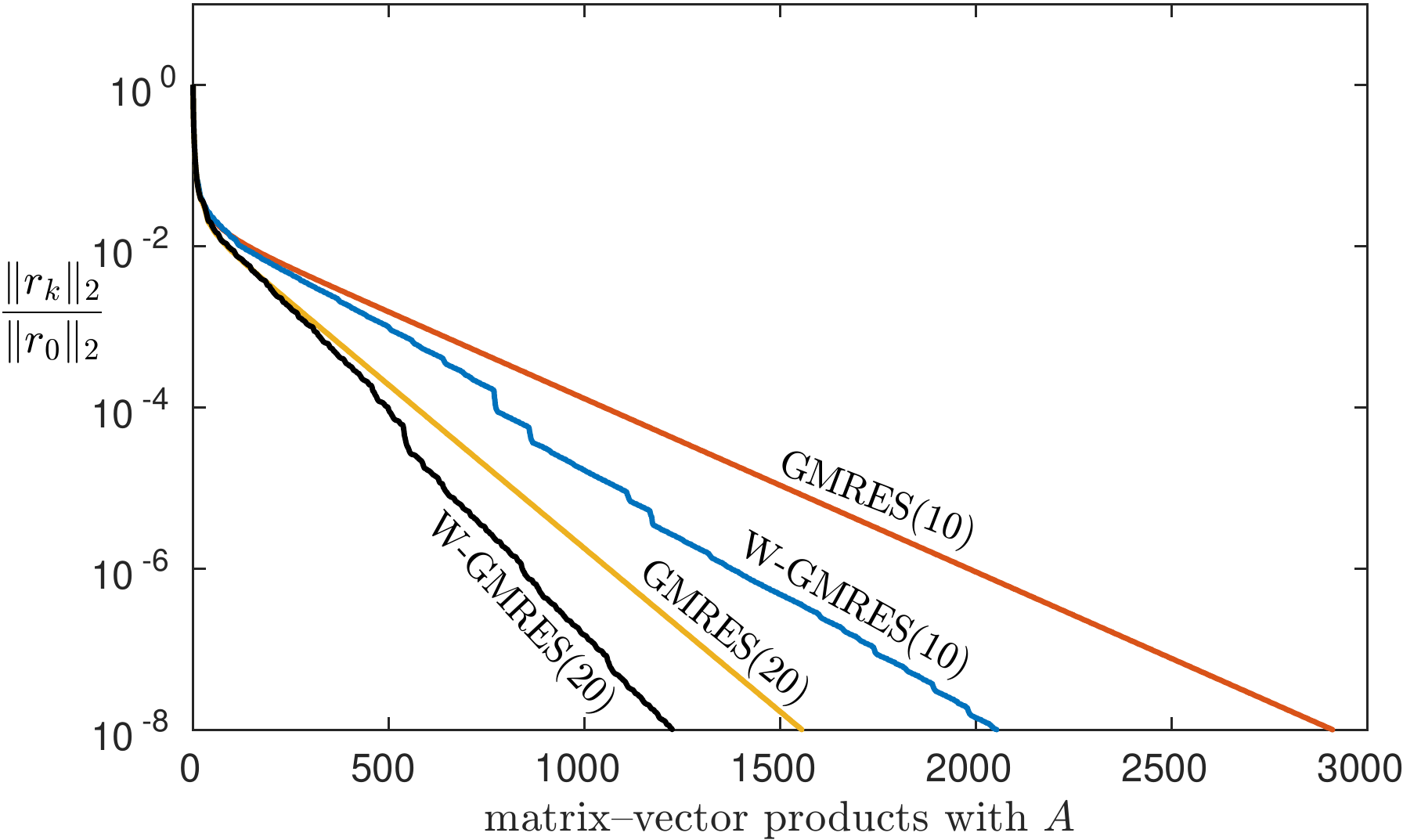} \end{center}

\vspace*{-4pt}
\caption{\label{fig:lap2d}
GMRES($m$) and W-GMRES($m$) for a discretization of the
Dirichlet Laplacian on the unit square.  The eigenvectors are not localized, 
but weighting still accelerates convergence.}
\end{figure}
\end{example}

\subsection{Extra weighting} \label{sec:extra}

If weighting can improve GMRES, especially when the smallest 
eigenvalues have localized eigenvectors, might more extreme weighting
help even more?
Recall from~\cref{eq:diag} that for diagonal $A$,
\begin{equation} \label{eq:diagextra}
 \|r_k\|_\W^2 = \mingmres \sum_{j=1}^n  |@\poly(\lambda_j)|^2 |s_j b_j|^2.
\end{equation}
Essai's standard weighting uses $s_j^2 = w_j = |b_j|$.
Suppose one instead takes $s_j^2 = w_j = |b_j|^p$ for $p\ge0$.
For large $p$ we find it particularly important to limit the minimal weight,
as in~\cref{eq:essai}; we require $w_j \in [10^{-10}, 1]$.

For most cases we have studied, extra weighting has little qualitative impact.  
For diagonal $A$, $p>1$ can accelerate convergence; 
for the general case, even with localized eigenvectors, results are mixed.
The next example shows a case where it helps.

\begin{example}[Extra weighting, localized eigenvectors] \label{ex:extra}\ \\
The \Orsirr\ matrix from Matrix Market~\cite{MM} is a nonsymmetric matrix of dimension $n=1030$.
\Cref{fig:comploc} shows the eigenvectors to be fairly well localized;
all the eigenvalues are real and negative.
\Cref{tbl:extra_orsirr} reports the matrix-vectors products (with $A$)
required by W-GMRES($m$) with different weighting powers~$p$
to satisfy the convergence criterion $\|r_k\|_2/\|r_0\|_2 \le 10^{-8}$.
(Here $r_0$ is a random Normal(0,1) vector.)
When $m$ is small (e.g., $m=5$ and $10$), extra weighting improves 
convergence significantly.
For example, when $m=10$, the extreme value $p=6$ gives the best result among
the reported~$p$.
For larger $m$, extra weighting has less impact.
\end{example}

\begin{table}[t!]
\caption{\label{tbl:extra_orsirr}
W-GMRES($m$) for the \Orsirr\ matrix, with extra weighting by 
the power $p$.
The table shows the matrix-vector products with $A$ required to reduce 
$\|r_k\|_2/\|r_0\|_2$ by a factor of $10^{-8}$.}

\begin{center}
\begin{tabular}{|r|r|r|r|r|r|r|r|r|r|r|}  \hline\hline
$p$      &     
\multicolumn{1}{|c|}{0}  &
\multicolumn{1}{|c|}{1}  &
\multicolumn{1}{|c|}{2}  &
\multicolumn{1}{|c|}{3}  &
\multicolumn{1}{|c|}{4}  &
\multicolumn{1}{|c|}{5}  &
\multicolumn{1}{|c|}{6}  &
\multicolumn{1}{|c|}{8}  &
\multicolumn{1}{|c|}{10}  \\ \hline \hline
$m= 5$   & 20000  & 20000  & 13854  &  8126  &  6641  &  7640  &  4043  &  3947  &  4512 \\ \hline
$m=10$   & 16299  & 13883  &  5009  &  5901  &  3810  &  3198  &  3053  &  3129  &  3064 \\ \hline
$m=20$   & 13653  &  2934  &  2569  &  2134  &  2221  &  2154  &  2303  &  2294  &  2377 \\ \hline
$m=30$   &  3750  &  2572  &  2224  &  2012  &  2075  &  2009  &  1998  &  1961  &  2094 \\ \hline
\end{tabular}
\end{center}
\end{table} 

\subsection{Random weighting} \label{sec:randw}

We have argued that 
when $A$ has local eigenvectors, residual-based weighting 
can encourage convergence in certain eigenvector directions
and break repetitive patterns in GMRES($m$).
In contrast, when the eigenvectors are predominantly global, weighting can still 
sometimes help by breaking patterns in GMRES($m$), but there is no
compelling connection to the eigenvectors.  
Thus if the eigenvectors of $A$ are global,
residual-based weighting should fare no better than the random 
diagonal weighting proposed by Saberi Najafi and Zareamoghaddam~\cite{SNZa}.  
In contrast, for localized eigenvectors, residual-based weighting should have a
significant advantage over random weighting.
The next two examples show this behavior.

\begin{example}[Random weights for a matrix with local eigenvectors]\ \\
The Add20 matrix from \cref{ex:add20} has localized eigenvectors, 
and W-GMRES($m$) improved convergence for $m=3$ and $m=6$.  
\Cref{tbl:random2} compares this performance to two random weighting schemes:
$W$ is a diagonal matrix with uniform random entries chosen in the interval $[0.5,1.5]$
(first method) and $[0,1]$ (second method), with fresh weights generated at each restart.
All experiments use the same $r_0$ as in \cref{fig:gwf1}; since performance changes 
for each run with random weights, we average the matrix-vector products with $A$ (iteration counts) over
ten trials each.  All three weighted methods beat GMRES($m$), but the residual-based
weighting has a significant advantage.
\end{example}

\begin{table}[t!]  \label{tbl:random2}
\caption{For the Add20 matrix, 
the matrix-vector products (with $A$) required to converge to
$\|r_k\|_2/\|r_0\|_2 \le 10^{-8}$ for
GMRES($m$), W-GMRES($m$), and 
random weights chosen uniformly in $[0.5,1.5]$ and $[0,1]$.
The random counts are averaged over 10~trials (same $r_0$).
W-GMRES($m$) does much better than the random schemes for this matrix with local eigenvectors.}

\begin{center}
\begin{tabular}{|r|r|r|r|r|r|}  \hline\hline
\multicolumn{1}{|c|}{$m$}  & 
\multicolumn{1}{|c|}{GMRES($m$)} & 
\multicolumn{1}{|c|}{W-GMRES($m$)}  & 
\multicolumn{1}{|c|}{rand(.5,1.5)}  & 
\multicolumn{1}{|c|}{rand(0,1)}   \\ \hline \hline
  1 & $>20000$ \hspace*{1em}&  1785\hspace*{2.25em}  & 15845.5\hspace*{.9em} & 4077.1 \hspace*{.25em}\\ 
  2 & $>20000$ \hspace*{1em} &  2882\hspace*{2.25em} & 10506.7\hspace*{.9em} & 8190.2 \hspace*{.25em}\\ 
  3 & $>20000$ \hspace*{1em} &  4580\hspace*{2.25em} & 8664.0\hspace*{.9em} & 7237.1 \hspace*{.25em}\\ 
  6 & 10032 \hspace*{1em}    &  1022\hspace*{2.25em} & 4270.7\hspace*{.9em} & 3382.7 \hspace*{.25em}\\ 
 10 &  4220 \hspace*{1em}    &   762\hspace*{2.25em} & 1945.4\hspace*{.9em} & 1331.9 \hspace*{.25em}\\ 
 15 &  2469 \hspace*{1em}    &   716\hspace*{2.25em} & 1531.6\hspace*{.9em} & 1126.8 \hspace*{.25em}\\ 
 20 &  1605 \hspace*{1em}    &   650\hspace*{2.25em} & 1109.1\hspace*{.9em} & 1022.5 \hspace*{.25em}\\   \hline\hline
\end{tabular}
\end{center}
\end{table}

\begin{table}[t!]  \label{tbl:random3}
\caption{For the discrete Laplacian on the unit square,
the matrix-vector products (with $A$) required to converge to
$\|r_k\|_2/\|r_0\|_2 \le 10^{-8}$ for
GMRES($m$), W-GMRES($m$), and two random schemes:
weights chosen uniformly in $[0.5,1.5]$ and $[0,1]$.
The random counts are averaged over 10~trials, all with the
same $r_0$.
For smaller values of $m$, the residual-based weighting of W-GMRES($m$) performs 
significantly worse than the random weighting for this matrix with global eigenvectors.}

\begin{center}
\begin{tabular}{|r|r|r|r|r|r|}  \hline\hline
\multicolumn{1}{|c|}{$m$}  & 
\multicolumn{1}{|c|}{GMRES($m$)} & 
\multicolumn{1}{|c|}{W-GMRES($m$)}  & 
\multicolumn{1}{|c|}{rand(.5,1.5)}  & 
\multicolumn{1}{|c|}{rand(0,1)}   \\ \hline \hline
  1 & $>20000$\hspace*{1em} & $>20000$\hspace*{1.5em} & 1446.0\hspace*{1em} & 1059.4\hspace*{0.5em} \\ 
  2 & 14449\hspace*{1em} &  9958\hspace*{1.5em} & 2499.8\hspace*{1em} & 2318.4\hspace*{0.5em} \\ 
  3 &  9687\hspace*{1em} &  6450\hspace*{1.5em} & 2670.6\hspace*{1em} & 2329.8\hspace*{0.5em} \\ 
  6 &  4867\hspace*{1em} &  3704\hspace*{1.5em} & 2267.1\hspace*{1em} & 1804.0\hspace*{0.5em} \\ 
 10 &  2912\hspace*{1em} &  2054\hspace*{1.5em} & 1726.3\hspace*{1em} & 1545.7\hspace*{0.5em} \\ 
 15 &  2017\hspace*{1em} &  1650\hspace*{1.5em} & 1742.0\hspace*{1em} & 1549.3\hspace*{0.5em} \\ 
 20 &  1556\hspace*{1em} &  1225\hspace*{1.5em} & 1404.1\hspace*{1em} & 1338.2\hspace*{0.5em} \\  \hline \hline
\end{tabular}
\end{center}
\end{table}

\begin{example}[Random weights for a matrix with global eigenvectors]\ \\
The results are very different for the discretization of the Dirichlet Laplacian 
on the unit square (\cref{ex:lap2d}).  
\cref{tbl:random3} shows that residual-based weighting
is inferior to most of the random weighting experiments, with particularly extreme
differences for small $m$.  
There appears to be no advantage to correlating weighting to residual components.
(This example shows a ``tortoise and hare'' effect~\cite{Em03}, 
where the random weighting with $m=1$ often converges in fewer iterations than
the larger $m$ here.)
\end{example}

\section{W-GMRES-DCT: Localizing eigenvectors with Fourier transforms} \label{sec:fft}

\begin{figure}[t!]
\begin{picture}(0,155)
 \put(0,0){\includegraphics[scale=0.5]{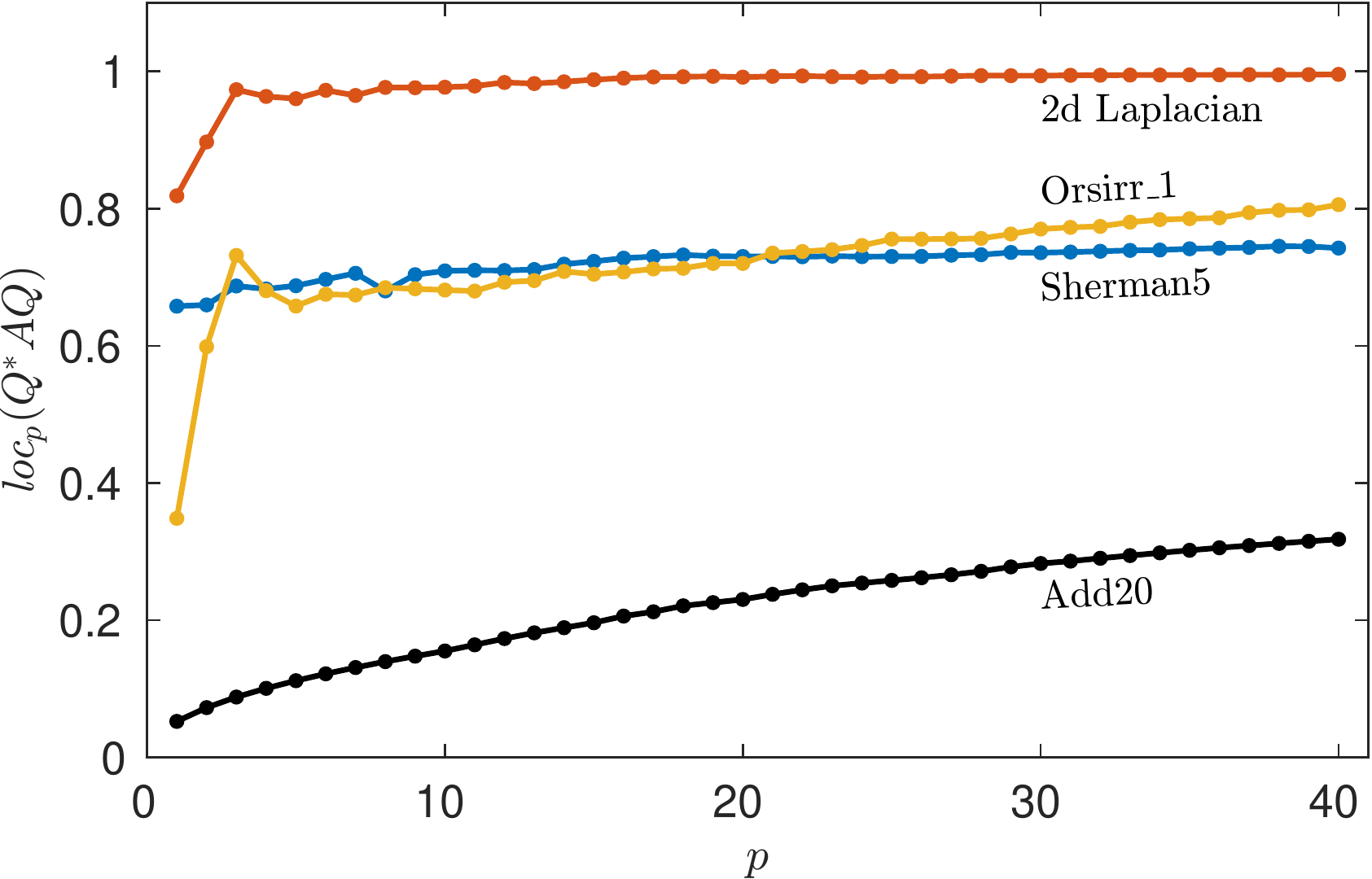}}
 \put(245,42){\footnotesize \begin{tabular}{lcc}
             \emph{matrix} &  $n$ & \hspace*{3.5em} \llap{$\loc_5(Q^*\!AQ)$} \\ \hline
              Add20 & $\!\!2395$ & $0.1118$ \\
              Orsirr\_1& $\!1030$ & $0.6579$\\
              Sherman5 & $\!\!3312$ & $0.6876$\\
              2d Laplacian\!\!\! & \!\!$9801$ & $0.9603$ \\
             \end{tabular}
            }
\end{picture}

\caption{\label{fig:comploc_dct}
Eigenvector localization measures ${\it loc}_p(Q^*AQ)$ for the 
coefficient matrices transformed by the DCT matrix $Q^* = {\tt dct(eye(n))}$.
The transformation reduces the localization of the Add20 eigenvectors,
while significantly enhancing localization for the three other matrices.
}
\end{figure}

Equation~\cref{eq:evweight} suggests that, for diagonalizable 
$A = V \Lambda V^{-1}$, the ideal weighting scheme should incorporate
eigenvector information.  It will be convenient to now express $W$ in the form
$W = Z^*Z$ with $Z = S V^{-1}$ for $ S = {\rm diag}(s_1, \ldots, s_n),$
giving
\[ \|r_k\|_\W^2 = \mingmres \sum_{j=1}^n  |@\poly(\lambda_j)|^2 |s_j b_j|^2.\]
Incorporating $V^{-1}$ into the weighting diagonalizes the problem and makes
residual weighting compelling, but such an inner product is clearly impractical.
However, we have seen that weighting can improve
GMRES convergence if the eigenvectors are merely localized.
This suggests a new weighting strategy that uses $W = Z^*Z$ with
\[ Z = S Q^*, \qquad S = {\rm diag}(s_1, \ldots, s_n),\]
where $Q$ is unitary 
and the weights $s_1, \ldots, s_n$ are derived from the residual as before.
Since GMRES on $(A,b)$ in this $W$-inner product amounts to standard 2-norm GMRES on
$(ZAZ^{-1},Zb) = (S Q^*\kern-2pt AQ S^{-1}, SQ^*b)$, we seek $Q$ in which 
\emph{the eigenvectors of $Q^*\kern-2pt AQ$ are localized}.
When $A$ is a discretization of a constant-coefficient differential operator, 
a natural choice for $Q^*$ is the discrete cosine transform matrix~\cite{Str99},%
\footnote{In Strang's parlance~\cite{Str99}, MATLAB's {\tt dct} 
routine gives a DCT-3, scaled so the matrix is unitary.}
since for many such problems the eigenvectors associated with the 
smallest-magnitude eigenvalues are \emph{low in frequency}.   
We call the resulting algorithm W-GMRES-DCT($m$).
In MATLAB, one could compute $Q^* = $ \verb|dct(eye(n))|,
a unitary but dense matrix that should never be formed; 
instead, it can be applied to vectors in $O(n\log n)$ time using the
Fast Fourier Transform.  In MATLAB code, the inner product becomes
\[ \langle x,y\rangle_W = y^*Q S^*S Q^* x = (SQ^*y)^*(SQ^*x) = \mbox{\texttt{(s.*dct(y))'*(s.*dct(x))}},\]
where ${\tt s} = [s_1, \ldots, s_n]^T$.
Alternatively, $Q^*$ could be the discrete Fourier transform matrix 
({\tt fft(eye(n))}, which uses complex arithmetic), 
a discrete wavelet transform,
or any other transformation that localizes (i.e., \emph{sparsifies})
the eigenvectors of $A$ (informed by the motivating application).
\Cref{fig:comploc_dct} shows how the DCT completely alters the
localization measure $\loc_p$~\cref{eq:locp} for the
matrices shown in \cref{fig:comploc}.

Presuming $Q$ is unitary, one can expedite the Arnoldi process 
by only storing the vectors $Q^*v_j$, rather than the usual Arnoldi vectors $v_j$.  
For example, in the $W$-inner product, the $k$th step of the 
conventional modified Gram--Schmidt Arnoldi process,
\medskip
\begin{center}\begin{minipage}{2.5in}
$\widehat{v} = A v_k$ \\ 
for $j=1,\ldots, k$\\
\hspace*{2em}$h_{j,k} = \langle \widehat{v}, v_j\rangle_W
                      =  v_j^*Q S^*\kern-1.5pt S Q^* \widehat{v} $ \\
\hspace*{2em}$ \widehat{v} = \widehat{v} - h_{j,k} v_j$  \\
end,
\end{minipage}
\end{center}
\medskip
can be sped up by storing only the vectors $u_j := Q^*v_j$, giving
\medskip
\begin{center}\begin{minipage}{2.5in}
$\widehat{u} = Q^* \kern-2pt A v_k = Q^*\kern-2pt  A Q@ u_k$ \\ 
for $j=1,\ldots, k$\\
\hspace*{2em}$h_{j,k} = \langle \widehat{v}, v_j\rangle_W
                      =  u_j^* S^*\kern-1.5pt S @\widehat{u}$ \\
\hspace*{2em}$ \widehat{u} = \widehat{u} - h_{j,k} u_j$  \\
end.
\end{minipage}
\end{center}
\medskip
This latter formulation, which we use in our experiments,
removes all applications of $Q$ from the innermost loop.
It is mathematically equivalent to the conventional version 
and differs only through unitary transformations,
but over many cycles one can observe small drift
in the numerical behavior of the two implementations.

\begin{figure}[t!]
  \begin{center}\includegraphics[scale=0.5]{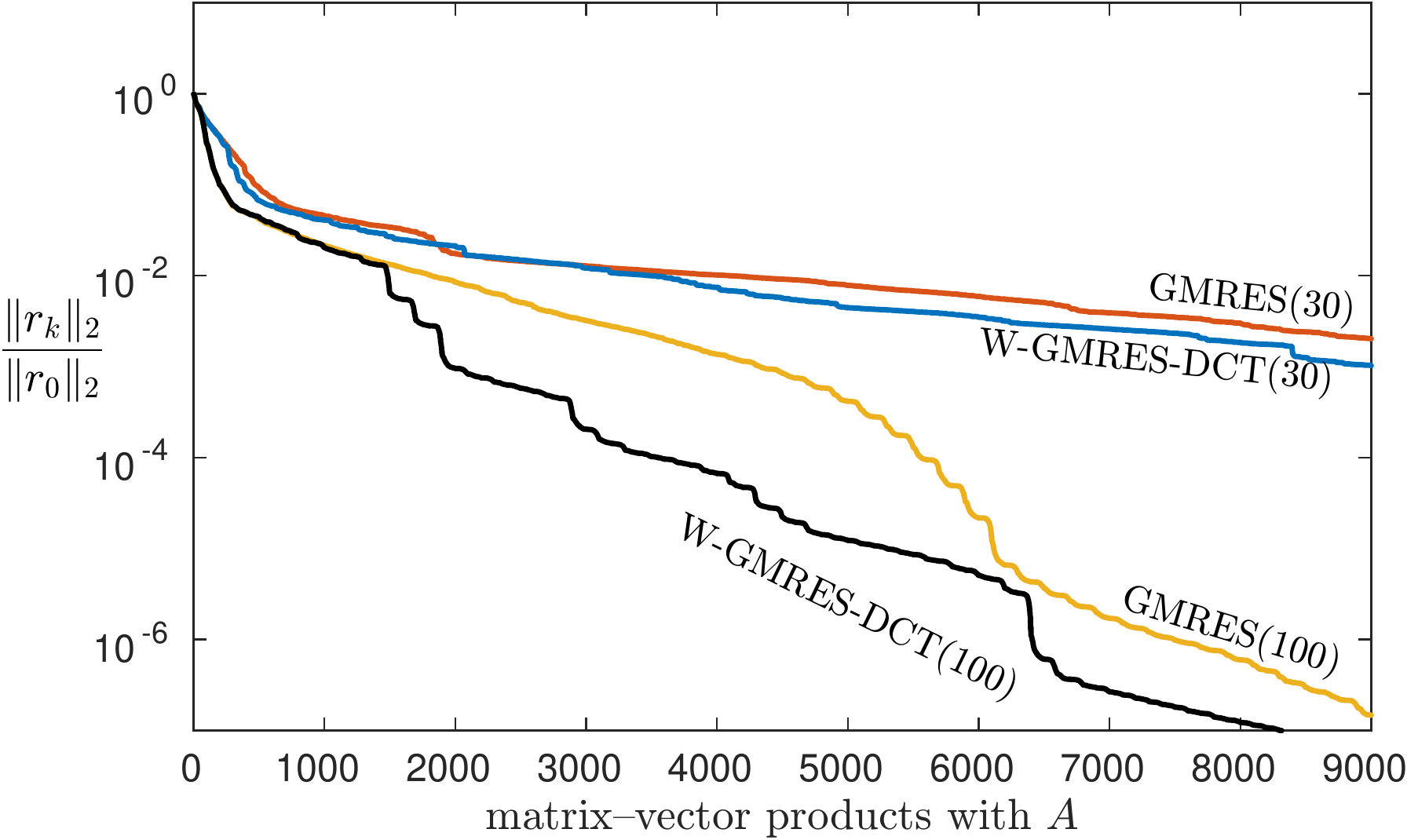} \end{center}

\vspace*{-4pt}
\caption{\label{fig:dct_sherex}
GMRES($m$) and W-GMRES-DCT($m$) for the nonsymmetric Sherman5 matrix.
The DCT helps to localize the eigenvectors, resulting in significantly
better convergence than seen for W-GMRES($m$) in the top plot 
of \cref{fig:sherex}.
}
\end{figure}

\begin{example}[Sherman5, revisited with discrete cosine transform]\ \\
How does the DCT affect the Sherman5 system (\cref{ex:sherman})?
Compare \cref{fig:dct_sherex} to the top of \cref{fig:sherex}; 
the DCT localizes the eigenvectors well, and now W-GMRES-DCT(30)
slightly outperforms GMRES(30), while W-GMRES-DCT(100) converges in 
fewer iterations than GMRES(100).  
The DCT does not localize the eigenvectors as effectively as 
diagonalization (bottom of \cref{fig:sherex}), but it is practical.
\end{example}


\begin{example}[Discrete Laplacian and convection--diffusion]\ \\
In \cref{ex:lap2d} we argued that W-GMRES($m$) converged 
better than GMRES($m$) despite global eigenvectors because
the weighted inner product disrupts the cyclic pattern of GMRES($m$).
The DCT inner product localizes the eigenvectors, and, as \cref{fig:dct_lap2d}
shows, convergence is improved.
Indeed, W-GMRES-DCT(20) requires roughly half as many iterations 
as W-GMRES(20), and both beat GMRES(20).

\begin{figure}[t!]
  \begin{center}\includegraphics[scale=0.5]{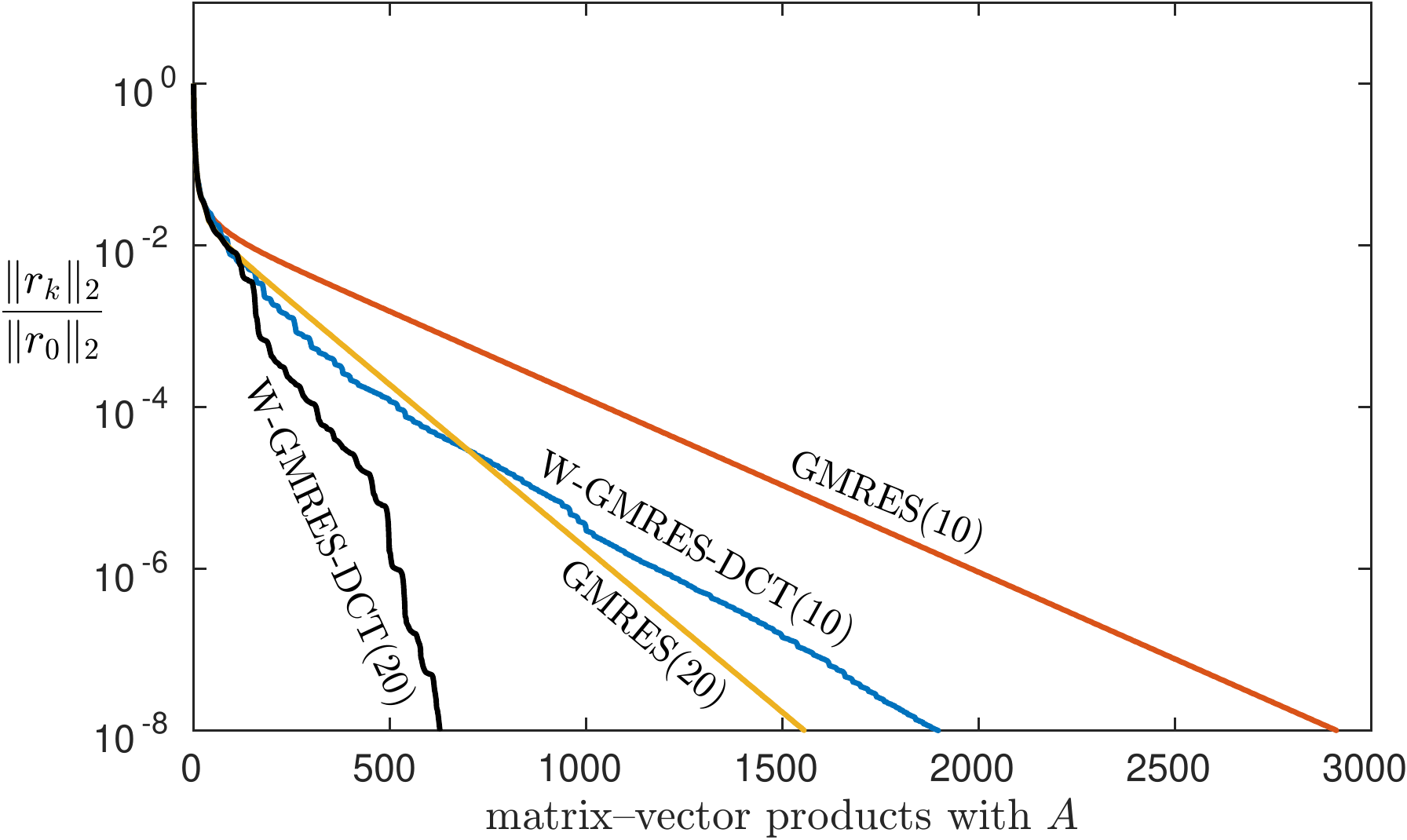} \end{center}

\vspace*{-4pt}
\caption{\label{fig:dct_lap2d}
GMRES($m$) and W-GMRES-DCT($m$) for the 2d Dirichlet
Laplacian (\cref{ex:lap2d}).
The DCT localizes the eigenvectors, improving the convergence
of W-GMRES($m$) in \cref{fig:lap2d}.
}
\end{figure}

The eigenvectors of the discrete Laplacian are two-dimensional sinusoids,
so the efficacy of the DCT is no surprise.  
The success of the method holds up when some nonsymmetry is added.
\Cref{fig:dct_convdiff} compares GMRES(10),
W-GMRES(10), and W-GMRES-DCT(10) for finite difference discretizations
of the convection--diffusion operator $L u = -(u_{xx}+u_{yy}) + u_x$,
again with Dirichlet boundary conditions and mesh size $h=1/100$ ($n=99^2=9801$).  
Due to the convection term $u_x$ the eigenvectors are no longer orthogonal, 
yet the DCT-enhanced weighting remains quite effective.

\begin{figure}[t!]
  \begin{center}\includegraphics[scale=0.5]{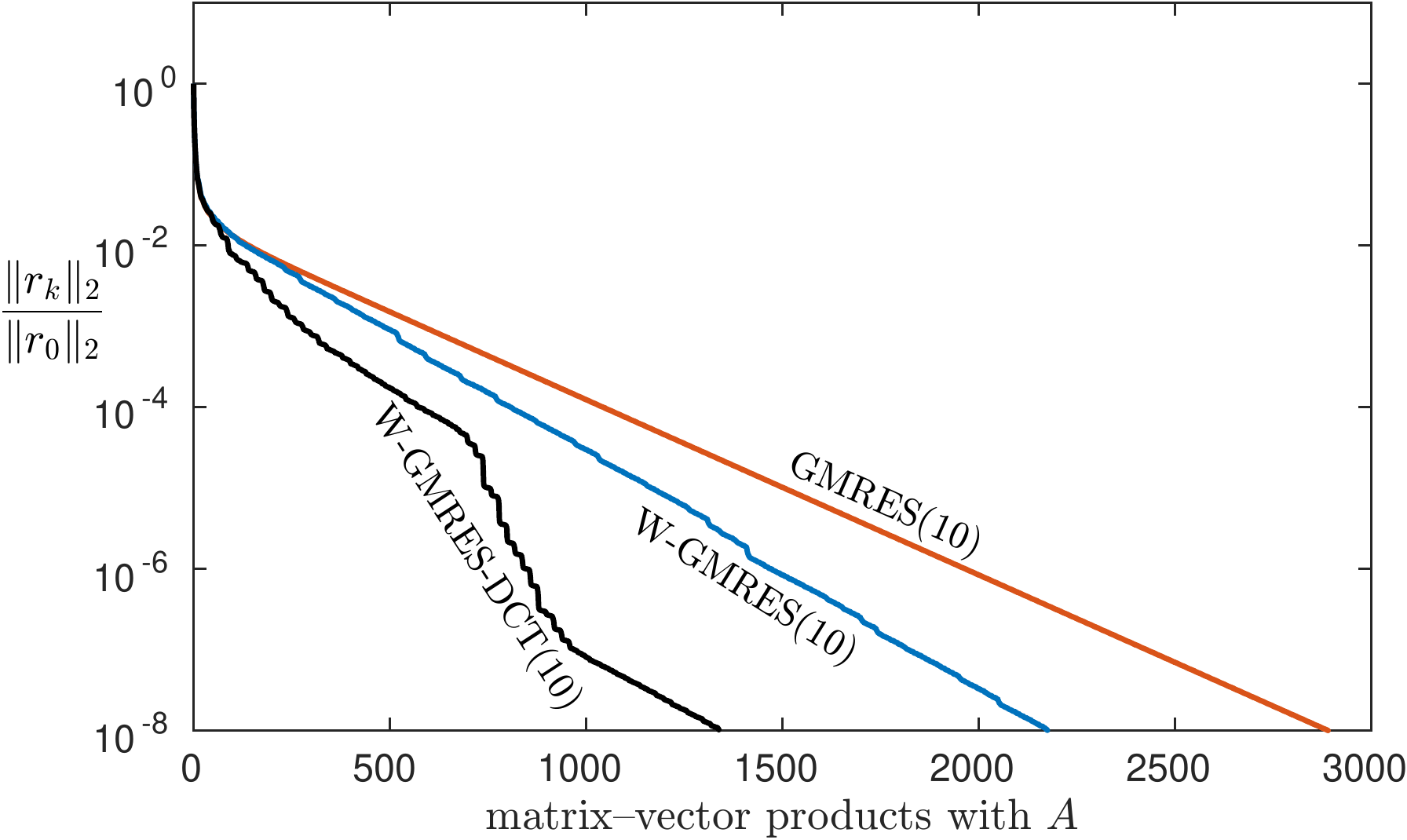} \end{center}

\vspace*{-4pt}
\caption{\label{fig:dct_convdiff}
GMRES(10), W-GMRES(10), and W-GMRES-DCT(10) for the finite difference
discretization of a 2d convection--diffusion problem with Dirichlet 
boundary conditions.
}
\end{figure}

\end{example}

\section{Weighted GMRES-DR} \label{sec:drwgmres}
To improve convergence, many modifications to restarted GMRES have been proposed;
see, e.g., \cite{BaJeKo,BaJeMa,DS99,Jo94b,Sa93,Si99,vdVVu}.  
Deflated GMRES methods use approximate eigenvectors in an attempt to
remove some troublesome eigenvalues from the spectrum;
see~\cite{Lan-DR,BaCaGoRe,BuEr,ChSa,ErBuPo,FrVu,GiDrPiVa,KhYe,GMRES-E,GMRES-DR,PadeStMaJoMa,Sa95B} 
and references in~\cite{Gu12}.  
Inner-product weighting has been incorporated into deflated GMRES in~\cite{Me09,NiLuZh}.  
Here we apply weighting to the efficient GMRES-DR method~\cite{GMRES-DR}.  

As detailed in~\cite{GMRES-DR}, GMRES-DR computes approximate eigenvectors 
while solving linear equations. 
We show how to restart with a change of inner product;
cf.~\cite[Prop.~1]{Es98}.  
The end of a GMRES-DR cycle in the $W$-inner product gives
\begin{equation}
 AV_k = V_{k+1} H_{k+1,k}, \label{eq:recur1}
\end{equation} 
where $V_k$ is an $n\times k$ matrix whose columns span a subspace of
approximate eigenvectors, $V_{k+1}$ is the same but with an extra column appended,
and $H_{k+1,k}$ is a full $(k+1)\times k$ matrix.  
Since this recurrence is generated with the $W$-inner product, the columns of
$V_{k+1}$ are $W$-orthogonal,
$V_{k+1}^* W V_{k+1}^{} = I$, and $H_{k+1,k}^{} = V_{k+1}^* W A V_{k}^{}$.  
To convert~\cref{eq:recur1} to a new inner product defined by $\widehat{W}$,  
compute the Cholesky factorization
\[
V_{k+1}^* \widehat{W} V_{k+1}^{} = R^*R,
\] 
where $R$ is a $(k+1)\times (k+1)$ upper triangular matrix. 
Let $R_k$ denote the leading $k\times k$ submatrix of $R$, and set
\[
\widehat{V}_{k+1} := V_{k+1} R^{-1}, \qquad
\widehat{H}_{k+1,k}^{} := R H_{k+1,k}^{} R_k^{-1}.\]
Then we have the new Arnoldi-like recurrence in the $\widehat{W}$-norm:
\[
 A \widehat{V}_k = \widehat{V}_{k+1} \widehat{H}_{k+1,k}
\] 
with $\widehat{V}_{k+1}^* \widehat{W} \widehat{V}_{k+1}^{} = I$ 
and $\widehat{H}_{k+1,k}^{} = \widehat{V}_{k+1}^* \widehat{W} A \widehat{V}_k^{}$.
Building on this decomposition, W-GMRES-DR extends the approximation subspace via the
Arnoldi process, and computes from this subspace the iterate that minimizes the
$\widehat{W}$-norm of the residual.

\begin{figure}[t!]
   \begin{center}\includegraphics[scale=0.5]{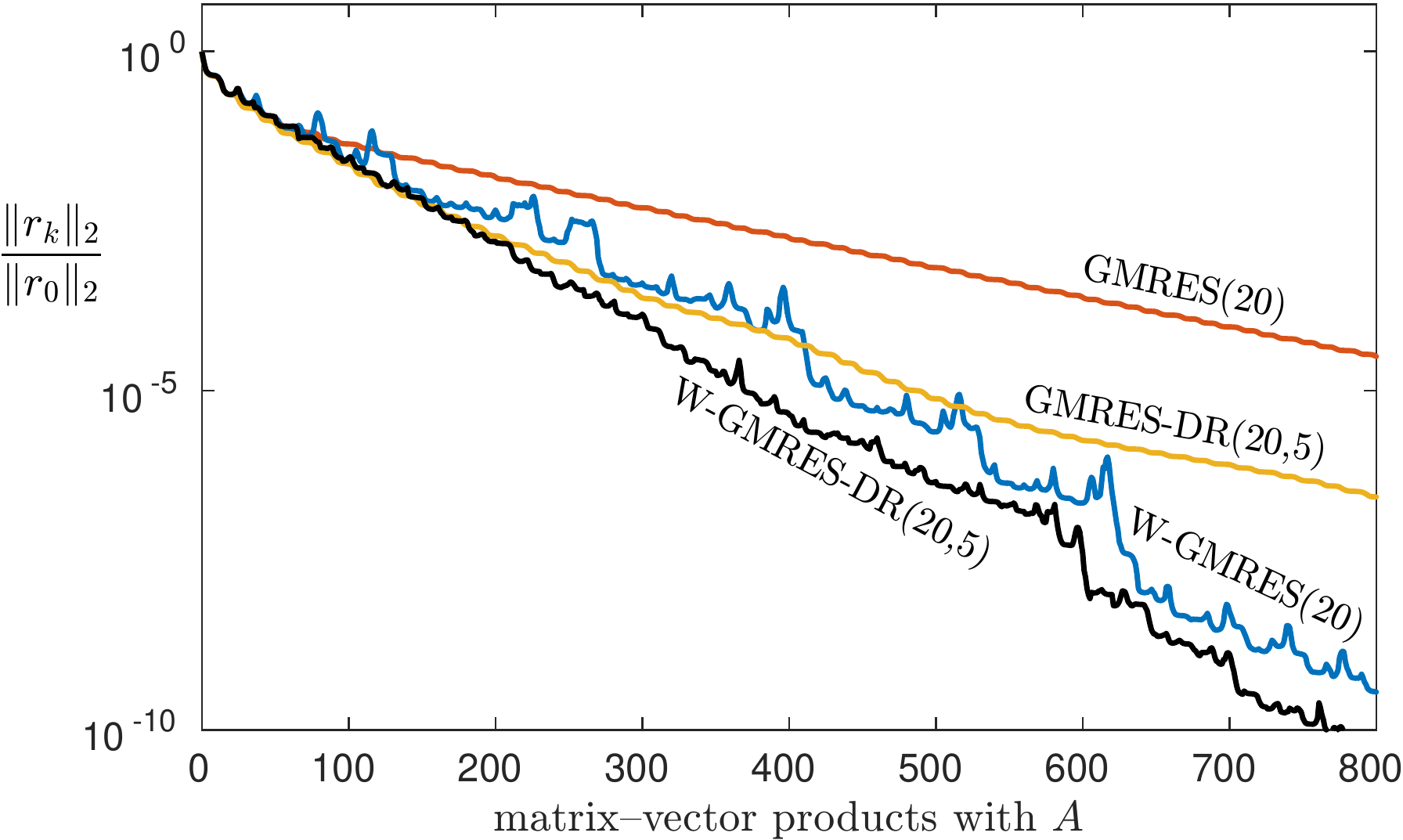} \end{center}

\vspace*{-5pt}
\caption{\label{fig:drex_add20}
Comparison of restarted GMRES, with and without weighting and deflation, for Add20.}
\end{figure}

\begin{example}
We apply W-GMRES-DR to the matrix Add20 from~\cref{ex:add20}.  
\Cref{fig:drex_add20} shows convergence for GMRES(20) and GMRES-DR(20,5), 
with and without weighting.  
Weighting improves convergence more than deflating eigenvalues;
combining the two strategies in W-GMRES-DR modestly improves convergence.  
\end{example}

\begin{figure}
\begin{center}\includegraphics[scale=0.5]{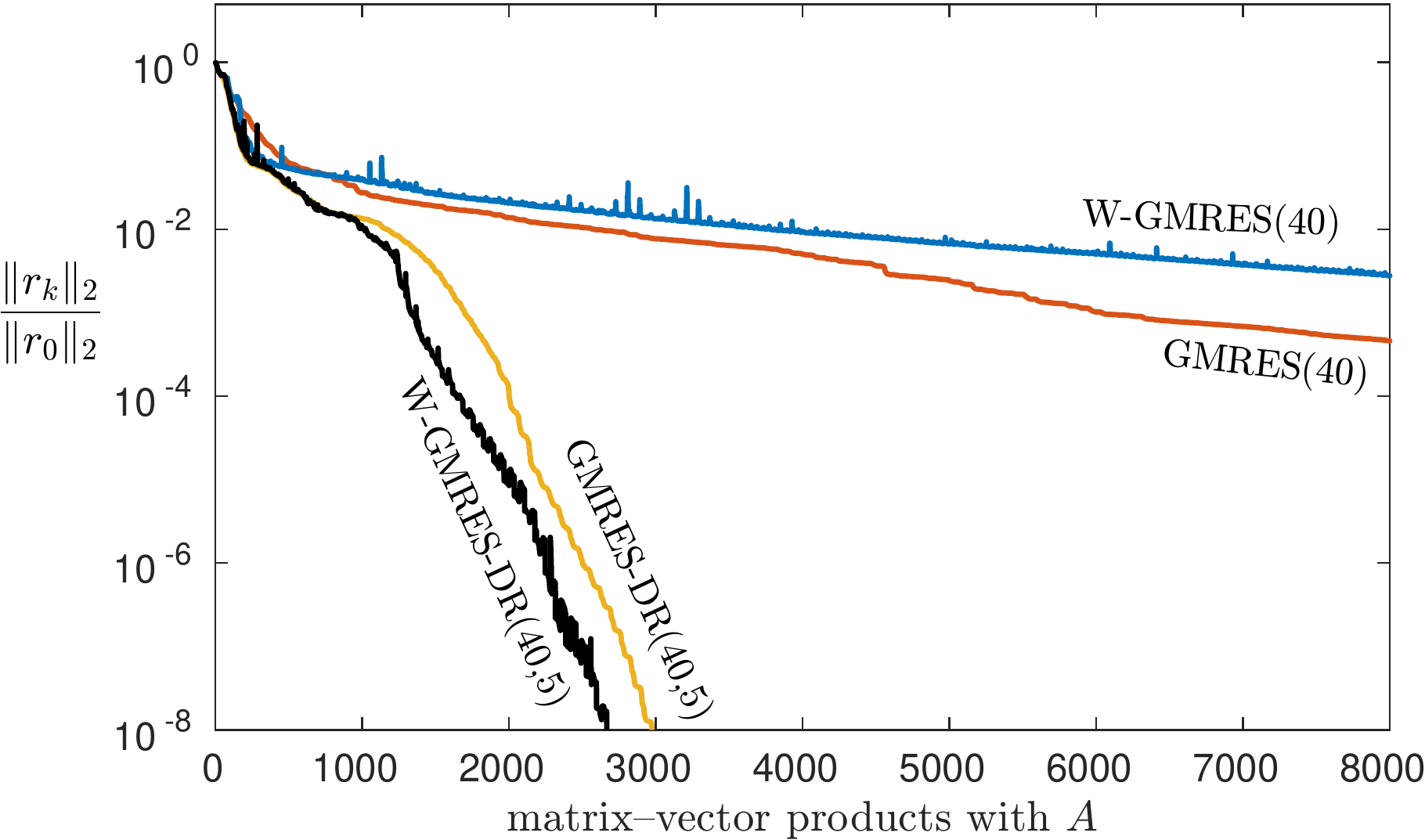} \end{center}

\vspace*{-5pt}
\caption{\label{fig:drex_sherex} 
Comparison of methods with and without weighting and deflation for Sherman5.}
\end{figure}

\begin{example}
For the matrix Sherman5, weighting did not help restarted GMRES in \cref{ex:sherman}.  
\Cref{fig:drex_sherex} adds GMRES-DR to the picture: deflating eigenvalues is 
very important, and now adding weighting to GMRES-DR(40,5) gives convergence in 
306~fewer matrix-vector products with $A$.
\end{example}

\section{Conclusions}

Since its introduction by Essai~\cite{Es98}, W-GMRES has
been an intriguing method whose adoption has been impeded by an 
incomplete understanding of its convergence properties.  
We have argued, through some simple analysis and 
numerous experiments,  that its advantage is due to (1)~its ability
to break repetitive cycles in restarted GMRES, and (2)~its ability 
to target small-magnitude eigenvalues that slow GMRES convergence,
\emph{provided the corresponding eigenvectors are localized}.
In cases where those vectors are not localized the proposed
W-GMRES-DCT algorithm, which combines residual weighting with a fast 
transform to sparsify eigenvectors, can potentially help.
Weighted inner products can also be applied to deflated GMRES methods, 
such as GMRES-DR.\ \  While the improvement appears to be less significant than
for standard GMRES, weighting can sometimes reduce the number of matrix-vector 
products even further than with the deflation alone or weighting alone.  

Many related topics merit further investigation,
such as additional sparsifying transformations or more general
non-diagonal weighting matrices, application of related ideas to 
Arnoldi eigenvalue computations (cf.~\cite{SG06}), and the use of weighting to 
avoid near-breakdown in the nonsymmetric Lanczos~\cite{NLan-DR} method.

\section*{Acknowledgments} 
We thank J\"org Liesen, Jennifer Pestana, Andy Wathen and two anonymous
referees for helpful comments about this work.
ME is grateful for the support of the Einstein Stiftung Berlin,
which allowed him to visit the Technical University of Berlin while working on this project.
RBM appreciates the support of the Baylor University Sabbatical Program.

\bibliographystyle{siam}
\bibliography{morgan}

\end{document}